\documentclass[11pt]{article}
\usepackage{amsmath}
\usepackage{amssymb}

\usepackage[german,english]{babel}
\usepackage{amsmath, latexsym}
\usepackage{amssymb}
\usepackage{url}
\usepackage{graphicx}
\usepackage[usenames]{color}
\usepackage{gastex}
\usepackage{longtable}
\usepackage{lscape}
\usepackage{tabularx}
\usepackage{multicol}
\usepackage{amsmath}
\usepackage{latexsym}
\usepackage{verbatim}
\usepackage{multirow}
\usepackage{ragged2e}
\newtheorem{lemma}{{\sc Lemma}}[section]
\newtheorem{theo}{{\sc Theorem}}[section]

\newtheorem{defn}{{\em Definition}}[section]

\newtheorem{conj}{{\sc Conjecture}}

\newtheorem{example}{{\sc Example}}

\begin{document}

\title{Hamiltonian Cycle in Semi-Equivelar Maps on the Torus}
\author{{Dipendu Maity  and Ashish Kumar Upadhyay}\\[2mm]
{\normalsize Department of Mathematics}\\{\normalsize Indian Institute of Technology Patna }\\{\normalsize Patliputra Colony, Patna 800\,013,  India.}\\
{\small \{dipendumaity, upadhyay\}@iitp.ac.in}}
\maketitle

\vspace{-5mm}

\hrule

\begin{abstract} Semi-Equivelar maps are generalizations of Archimedean solids to the surfaces other than 2-sphere. There are eight semi-equivelar maps of types $\{3^{3},4^{2}\}$, $\{3^{2},4,3,4\}$, $\{6,3,6,3\}$, $\{3^{4},6\}$, $\{4,8^{2}\}$, $\{3,12^{2}\}$, $\{4,6,12\}$, $\{6,4,3,4\}$ exist on the torus. In this article we show the existence of Hamiltonian cycle in each semi-equivelar map on the torus except the map of type $\{3,12^{2}\}$. This result gives the partial solution to the conjecture which is given by Gr$\ddot{u}$nbaum \cite{grunbaum} and Nash-Williams \cite{nash williams} that every 4-connected graph on the torus is Hamiltonian.
\end{abstract}

{\small

{\bf AMS classification\,: 52B70, 05C45, 52C38}

{\bf Keywords\,:} Semi-Equivelar Maps, Cycles , Hamiltonian Cycles.
}

\bigskip

\hrule

\section{Introduction and definitions}

Let $X$ and $Y$ be two finite abstract simplicial complexes. An {\em isomorphism} from $X$ to $Y$ is a bijection $\phi : V(X)\rightarrow V(Y)$ such that for $\sigma \subseteq V(X)$, $\sigma$ is a simplex of $X$ if and only if $\phi(\sigma)$ is a simplex of $Y$. Two simplicial complexes $X$, $Y$ are called (simplicially) isomorphic when such an isomorphism exists. An isomorphism from a simplicial complex $X$ to itself is called an $automorphism$ of $X$. All the automorphisms of $X$ form a group, which is denoted by $Aut(X)$. A simplicial complex $M$ is said to be a {\em triangulated d-manifold} if it's geometric carrier $|M|$ is a topological $d$-manifold. When $d =2$, a connected combinatorial $2$-manifold $X$ is said to be vertex-transitive if the automorphism group of $X$ acts transitively on $V(X)$.

Let $G_{1}(V_{1}, E_{1})$ and $G_{2}(V_{2}, E_{2})$ be two simple graphs on the torus. Then $G_{1}\cup G_{2}$ is a graph $G(V, E)$ where $V = V_{1}\cup V_{2}$ and $E = E_{1}\cup E_{2}$. Similarly, $G_{1}\cap G_{2}$ is a graph $G(V, E)$ where $V = V_{1}\cap V_{2}$ and $E = E_{1}\cap E_{2}$. A {\em map} $M$ is an embedding of a graph $G$ on a surface $S$ such that the closure of components of $S \setminus G$, called the faces of $M$, are closed $2$-cells i.e. each homeomorphic to $2$-disk. A map $M$ is said to be a polyhedral map if the intersection of any two distinct faces is either empty, a common vertex, or a common edge. The following definitions are given more clearly in \cite{upadhyay1}. A {\em a-cycle} $C_{a}$ is a finite connected {\em 2-regular} graph with {\em p} vertices and the {\em face sequence} of a vertex {\em v} in a map is a finite sequence $(a^{p},b^{q},\dots,m^{r})$ of powers of positive integers $a, b, \dots, m \geq 3$ and $p, q, \dots, r \geq 1$ such that through the vertex {\em v, p} number of $C_{a}$ ($C_{a}$ denote the $a-cycle$), {\em q} number of $C_{b}$, \dots, {\em r} number of $C_{m}$ are incident. A map {\em K} is said to be {\em semi-equivelar} if face sequence of each vertex is same. Let {\em K} be a map of type $(a^{p},b^{q},\dots,m^{r})$ on $n$ vertices on the torus. By {\em Euler's} formula, $n - \frac{nd}{2}+(\frac{nd_{p}}{p}+\frac{nd_{q}}{q}+\dots+\frac{nd_{r}}{r}) = 0 $, that is, $1 - \frac{d}{2}+(\frac{d_{p}}{p}+\frac{d_{q}}{q}+\dots+\frac{d_{r}}{r}) = 0$, that is, $\frac{d}{2}-1 = \frac{d_{p}}{p}+\frac{d_{q}}{q}+\dots+\frac{d_{r}}{r}$, where {\em d} denote the degree of a vertex and $d_{t}$ denote the number of adjacent faces of length $t$ at a vertex for $t \in \{ p, q,\dots, r\}$. Now for fixed $d$, we check all possible values of $d_{t}$ for $t \in \{ p, q,\dots, r\}$ such that $\frac{d}{2}-1 = \frac{d_{p}}{p}+\frac{d_{q}}{q}+\dots+\frac{d_{r}}{r}$. This gives eight semi-equivelar maps of types $\{3^{3},4^{2}\}$, $\{3^{2},4,3,4\}$, $\{6,3,6,3\}$, $\{3^{4},6\}$, $\{4,8^{2}\}$, $\{3,12^{2}\}$, $\{4,6,12\}$, $\{6,4,3,4\}$. These maps are vertex-transitive by the result of \cite{such}. This imply, the edge graph of each semi-equivelar map on the torus is vertex-transitive. The connectivity number $k(G)$ of a graph $G$ is defined as the minimum number of vertices whose removal from $G$ results in a disconnected graph or in the trivial graph. A graph $G$ is said to be $k$-connected if $k(G)$ $\geq$ $k$. Now, we introduce the result given in \cite{godsil_royle}, a vertex-transitive graph with degree $k$ has connectivity at least $\frac{2}{3}(k+1)$. This gives the maps of types $\{4,8^{2}\}$, $\{3,12^{2}\}$, $\{4,6,12\}$ are $3$-connected, $\{6,3,6,3\}$, $\{6,4,3,4\}$ are $4$-connected and $\{3^{3},4^{2}\}$, $\{3^{2},4,3,4\}$, $\{3^{4},6\}$ are at least $4$-connected. Therefore the maps of types $\{3^{3},4^{2}\}$, $\{3^{2},4,3,4\}$, $\{6,3,6,3\}$, $\{6,4,3,4\}$, $\{3^{4},6\}$ are $4$-connected. \noindent Gr$\ddot{u}$nbaum \cite{grunbaum} and Nash-Williams \cite{nash williams} conjectured the following:

\begin{conj}\label{conj1} Every $4$-connected graph on the torus has a Hamiltonian cycle.
\end{conj}

We summarize the known existence results about existence of Hamiltonian cycles in graphs on torus and Klein bottle these are the some partial results of the conjecture \ref{conj1}. Altshuler \cite{a1} has shown that every $6$ and $4$-connected equivelar maps of type $\{3, 6\}$ and $\{4, 4\}$ respectively on torus have Hamiltonian cycles. Burnet and Richter\cite{brunet richter} have shown that every $5$-connected  triangulations on torus is Hamiltonian and Thomas and Yu\cite{thomas yu} have improved this result for any $5$-connected graphs. Brunet, Nakamoto and Negami \cite{brunet nakamoto negami} have shown every $5$-connected triangulated Klein bottle is Hamiltonian. Kawarabayashi and Ozeki\cite{kawarabayashi ozeki} showed every $4$-connected triangulated torus is Hamiltonian. In this article, we give partial solution of the conjecture as follows,\

\begin{theo}\label{thm1} Every $4$-connected semi-equivelar maps on torus has a Hamiltonian cycle.
\end{theo}

\noindent{\sc Proof of Theorem}\ref{thm1}: Proof of the theorem\ref{thm1} follows from theorem \ref{thm1.3}, \ref{thm2.1}, \ref{thm3.2}, \ref{thm4.1}, \ref{df8.1}.
\hfill$\Box$

Main result of this article is the following theorem.

\begin{theo}\label{thm2} The map of types $\{3^{3},4^{2}\}$, $\{3^{2},4,3,4\}$, $\{6,3,6,3\}$, $\{6,4,3,4\}$, $\{3^{4},6\}$, $\{4,8^{2}\}$, $\{4,6,12\}$ on the torus is Hamiltonian.
\end{theo}

\noindent{\sc Proof of Theorem}\ref{thm2}: Proof of the theorem \ref{thm2} follows from theorem \ref{thm1}, \ref{df5.1}, \ref{thm7.1}.
\hfill$\Box$

\begin{theo}\label{thm3} The map of type $\{3, 12, 12\}$ on the torus is non-Hamiltonian.
\end{theo}

\noindent{\sc Proof of Theorem}\ref{thm3}: Proof of the theorem \ref{thm3} follows from theorem \ref{thm6.1}.
\hfill$\Box$\\

Let {\em C, D} be two cycles. We define a {\em cylinder} $S_{C,D}$ be a set of complex with two boundary components {\em C, D}. This notion has been introduced in \cite{upadhyay1}. Let {\em C, D} be two cycles of same type in a map on the torus. The type of the cycles has been defined clearly in each section. The cycle $C$ is said to be {\em homologous} to $D$ if there exist a cylinder which is bounded by $C$ and $D$. We denote a cycle $u_{1}\rightarrow u_{2}\rightarrow \dots \rightarrow u_{k}\rightarrow u_{1}$  by $C(u_{1},u_{2},\dots,u_{k})$ in this article. Let $C_{1}(u_{1},u_{2},\dots,u_{r})$ and $C_{2}(v_{1},v_{2},\dots,v_{t})$ be two cycles. Let $u_{i}u_{i+1} = v_{k}v_{k+1}$. We define a cycle $C(u_{1},$ $\dots,$ $u_{i},$ $v_{k-1},$ $v_{k-2},$ $\dots,$ $v_{1},$ $v_{t},$ $v_{t-1},$ $\dots,$ $v_{k+1},$ $u_{i+2},$ $\dots,$ $u_{r})$. If $C_{1}$ and $C_{2}$ both bounds $2$-disk then $C$ bounds a $2$-disk. The cycle $C$ is called {\em concatenation} of $C_{1}$ and $C_{2}$. Let $Q_{1}(u_{1}\rightarrow u_{k})$ be a path. We call a path $Q_{2}(v_{1}\rightarrow v_{r})$ is extended from $Q_{1}$ if $V(Q_{1})\subset V(Q_{2}),$ $E(Q_{1})\subset E(Q_{2})$ and $v_{1}, v_{r} \not \in V(Q_{1})$. We call $Q_{2}$ is {\em extended path}. A notion of a {\em normal path} and {\em normal cycle} are defined in \cite{a1} to show the existence of Hamiltonian cycle in $4$ and $6$-regular maps on the torus. We have used the analogous definitions of {\em normal path} and {\em normal cycle} in this article. We define paths and cycles to discuss the existence of Hamiltonian cycles in each map of types $\{3^{3},4^{2}\}$, $\{3^{2},4,3,4\}$, $\{6,3,6,3\}$, $\{3^{4},6\}$, $\{4,8^{2}\}$, $\{3,12^{2}\}$, $\{4,6,12\}$, $\{6,4,3,4\}$ on the torus. The basic idea introduced as follows $\colon$

Let $M$ be a semi-equivelar map on the torus. We cut the torus along meridian and longitude, we get a planar representation. We define cycles in $M$ at a vertex locally such that each divides the degree of the vertex into a fixed ratio and maintains this ratio along the sequence of vertices throughout the cycle. We cut $M$ along any two non homologous cycles at a vertex to get a planer representation of $M$. We denote this representation by {\em $(r, s, k)$-representation}. In this article we use $T(r,s,k)$ to represent $(r, s, k)$-representation of a map. This representation has been defined more clearly in each section. This representation is called a planar representation of the map $M$ on the torus. This $(r, s, k)$-representation always exist for a semi-equivelar map. This idea also has been used in \cite{maity upadhyay1}. A cycle which contains all the vertices of the map and bounds a $2$-disk is called {\em contractible Hamiltonian} cycle. If it does not bounds any $2$-disk then we call {\em non-contractible}. We show the existence of contractible and non-contractible Hamiltonian cycles using $T(r, s, k)$ of $M$. The concept of contractible Hamiltonian cycles also has been studied in \cite{maity upadhyay0, upadhyay0}. We use this contractible Hamiltonian cycle to show the existence of Hamiltonian cycle in $M$.

\section{Map of type $\{3^{3}, 4^{2}\}$}

Let $M$ be a map of type $\{3^{3}, 4^{2}\}$ on the torus. We define paths as follows $\colon$

\begin{defn}\label{def1.1} Let $P(u_{1}\rightarrow\dots\rightarrow u_{i-1}\rightarrow u_{i}\rightarrow u_{i+1}\rightarrow\dots \rightarrow u_{r})$ be a path in $M$. The path $P$ is said to be well defined at an inner vertex $u_{i}$ if adjacent all the triangles on one side and all quadrangles on the other side of $u_{i-1}\rightarrow u_{i}\rightarrow u_{i+1}$ at $u_{i}$. If $u_{t}$ is a boundary vertex of $P$ then there is a well defined extended path of $P$ where $u_{t}$ is an inner vertex. We denote this path by $A_{1}$.
\end{defn}

We define link of a vertex $w$ and denote it by $lk(w)$. If $lk(w)$ contains any bold vertex then it indicates the non adjacent vertex to $w$. In this section, we consider a vertex $w$ with link $lk(w)=C(\textbf{v}_{2},v_{3},\textbf{v}_{4},v_{5},v_{6},v_{7},v_{1})$. In this link $\textbf{v}_{2}$ and $\textbf{v}_{4}$ are non adjacent vertices. We consider the permutation of vertices in the link $lk(w)$ counter clockwise locally at $w$. We use the notation throughout this article.

\begin{defn}\label{def1.2}
Let $P(u_{1}\rightarrow\dots\rightarrow u_{i-1}\rightarrow u_{i}\rightarrow u_{i+1}\rightarrow\dots \rightarrow u_{r})$ be a path in $M$. Let $u_{i}, u_{i+1}$ be two inner vertices of $P$ or an extended path of $P$. Then $P$ is said to be well defined at vertex $u_{i}$ if
\begin{enumerate}
 \item  the link $lk(u_{i})=C(\textbf{a},$ $u_{i-1},$ $\textbf{b},$ $c,$ $u_{i+1},$ $d,$ $e)$ gives the link $lk(u_{i+1})=C(\textbf{a}^{'},$ $u_{i+2},$ $\textbf{b}^{'},$ $d,$ $u_{i},$ $c,$ $p)$ and
 \item the link $lk(u_{i})=C(\textbf{x},u_{i+1},\textbf{z}, l, u_{i-1}, k, m)$ gives the link $lk(u_{i+1})=C(\textbf{x}^{'},u_{i},$ $\textbf{z}^{'},$ $k,$ $u_{i+2},$ $l,$ $p)$.

 \end{enumerate}
 We called $u_{i+1}$, $u_{i-1}$ are successor and predecessor of $u_{i}$ respectively in $P$. We denote this path by $A_{2}$.
\end{defn}

\begin{defn}\label{def1.3}
Let $P(u_{1}\rightarrow\dots\rightarrow u_{i-1}\rightarrow u_{i}\rightarrow u_{i+1}\rightarrow\dots \rightarrow u_{r})$ be a path in $M$. Let $u_{i}, u_{i+1}$ be two inner vertices of $P$ or an extended path of $P$. Then $P$ is said to be well defined at vertex $u_{i}$ if
\begin{enumerate}
 \item  the link $lk(u_{i})=C(\textbf{a},$ $u_{i-1},$ $\textbf{b},$ $c,$ $d,$ $u_{i+1},$ $e)$ gives the link $lk(u_{i+1})=C(\textbf{a}^{'},$ $u_{i+2},$ $\textbf{b}^{'},$ $p,$ $e,$ $u_{i},$ $d)$ and
 \item the link $lk(u_{i})=C(\textbf{a}^{'},$ $u_{i+1},$ $\textbf{b}^{'},$ $p,$ $e,$ $u_{i-1},$ $d)$ gives the link $lk(u_{i+1})=C(\textbf{a},$ $u_{i},$ $\textbf{b},$ $c,$ $d,$ $u_{i+2},$ $e)$.

 \end{enumerate}
 We denote this path by $A_{3}$.
\end{defn}

In each section, we have defined paths in map of fixed type. We have used the different notation for different path at each definition. In this section, we have consider the paths of type $A_{t}$ for $t \in \{1, 2, 3\}$. Now we consider maximal path of type $A_{t}$ for $t \in \{1, 2, 3\}$. We show that they are cycles. More precisely, let $Q(v_{1}\rightarrow v_{r})$ be a maximal path. Then, by the property of maximality, there is an edge $v_{r}v_{1}$ in the map. This gives the cycle $C(v_{1}, v_{2}, \dots, v_{r})$. Therefore we claim that,\

\begin{theo}\label{thm1.1}  Maximal path of type $A_{t}$ is cycle for $t \in \{1, 2, 3\}$.
\end{theo}
\noindent{\sc Proof of Theorem}\ref{thm1.1}  Let $P(u_{1}\rightarrow u_{2}\rightarrow ...\rightarrow u_{r})$ be a maximal path of type $A_{1}$. We consider the vertex $u_{r}$ and link $lk(u_{r})=C( \textbf{x},y,\textbf{z},w,v,u,u_{r-1})$. If $w = u_{1}$ then $C(u_{1}, u_{2}, \dots, u_{r})$ is a cycle. If $w \neq u_{1}$ then by the definition \ref{def1.1} $u_{r}$ is an inner vertex in the extended path of $P$. This gives a path, say $Q$, which is extended from $P$ and $length(P) < length(Q)$. Which is a contradiction as $P$ is  maximal path. Therefore $w = u_{1}$ and the path $P$ gives the cycle $C(u_{1}, u_{2}, \dots, u_{r})$. We use similar argument to show every maximal path of types $A_{2}$ and $A_{3}$ is cycle.
\hfill$\Box$\\

In each section, we use the terminology cycle in place of maximal path which is cycle. In this section it is cycle by theorem \ref{thm1.1}. Now, we show that if any two cycles of type $A_{t}$ intersect then they are identical for $t \in \{1, 2, 3\}$. Therefore we claim that,\

\begin{lemma}\label{lem1.1} Let $C_{1}, C_{2}$ be two cycles of type $A_{t}$ for $t \in \{1, 2, 3\}$. Then, $C_{1} = C_{2}$ if $C_{1} \cap C_{2} \neq \emptyset$.
\end{lemma}

\noindent{\sc Proof of Lemma}\ref{lem1.1} Let $C_{1}(u_{1_{1}},u_{1_{2}},\dots,u_{1_{r}})$ and $C_{2}(u_{2_{1}}, u_{2_{2}},\dots, u_{2_{s}})$ be two cycles of type $A_{1}$. If $C_{1} \cap C_{1} \neq \emptyset$ then there is a vertex $w \in V(C_{1} \cap C_{2})$. The cycles $C_{1}$ and $C_{2}$ both are well defined at $w$. We consider the link $lk(w) = C(\textbf{w}_{1}, w_{2}, \textbf{w}_{3},w_{4},w_{5},w_{6},w_{7})$. By definition \ref{def1.1}, $w_{4}, w_{7} \in C_{t}$ for $t \in \{1, 2\}$. This gives $w_{4}\rightarrow w\rightarrow w_{7}$ is part of $C_{t}$ for $t \in \{1, 2\}$. If $w = u_{1_{t_{1}}} = u_{2_{t_{2}}}$ then $w_{4} = u_{1_{t_{1}-1}} = u_{2_{t_{2}-1}}$ and $w_{7} = u_{1_{t_{1}+1}} = u_{2_{t_{2}+1}}$ for some $t_{1} \in \{1,\dots,r\}$ and $t_{2}\in \{1,\dots,s\}$. Now, we have $u_{1_{t_{1}-1}} = u_{2_{t_{2}-1}}$. We continue with this argument. We stop after finite steps as $r, s$ both finite. If $ r < s $ then we get $u_{1_{1}} = u_{2_{I+1}}, u_{1_{2}} = u_{2_{I+2}},\dots, u_{1_{r}} = u_{2_{I+r}}$ and $u_{1_{1}} = u_{2_{I+r+1}}$ for some $I \in \{1,\dots,s\}$. This imply that $I+1 = I+r+1$ and the cycle $C_{2}$ which contains a cycle of length $r$. This gives $r = s$. Hence $C_{1} = C_{2}$. We use similar argument for the cycles of type $A_{t}$ for $t \in \{2,3\}$.
\hfill$\Box$\\

We choose a cycle $C$ of type $A_{t}$ for $t \in \{1, 2, 3\}$ in $M$. We define a cylinder $S_{C}$ as follows. Let $u$ be a vertex in $C$. We collect all adjacent faces of $u$ in a set, say $S$.  We do it for each vertex of $C$. We consider geometric carrier of $S$ which is a cylinder and denote it by $S_{C}$. It has two boundary cycles either disjoint or identical by lemma \ref{lem1.1}. Let $\partial S_{C} = \{C_{1}, C_{2}\}$. We claim that $length(C) = length(C_{1}) = length(C_{2})$, that is,

\begin{lemma} \label{lem1.2} Let $S_{C}$ be a cylinder where $C$ is a cycle of type $A_{t}$ for $t \in \{1,2,3\}$. Also, let $\partial S_{C} = \{C_{1}, C_{2}\}$. Then $length(C) = length(C_{1}) = length(C_{2})$.
\end{lemma}

\noindent{\sc Proof of Lemma}\ref{lem1.2} Let $S_{C}$ be a cylinder of type $A_{1}$. Let $F_{1}, F_{2}, \dots, F_{r}$ be a sequence of adjacent faces in one side of $C$. The same sequence $F_{1}, F_{2}, \dots, F_{r}$ of faces is also adjacent to $C_{t}$ on one side of $C_{t}$ for $ t = 1$ or $2$. Therefore the sequence of types of faces in the other side of $C$ will follow on the other side of $C_{t}$ for $ t = 1$ or $2$. This gives the cycle $C_{t}$ is of type $A_{1}$ for $t \in \{1, 2\}$. It can be seen clearly from Figure 1 by considering the cycle $C(x_{1},\dots,x_{r})$. Here $S_{C}=\{C(w_{1},\dots,w_{r}), C(z_{1},\dots, z_{r})\}$. Similar argument we use for the cycle of type $A_{j}$ for $j \in \{2, 3\}$.

Let $length(C) \neq length(C_{1}) \neq length(C_{2})$. Also, let $ C(u_{1}, u_{2},\dots, u_{r})$, $ C_{1}(v_{1},$ $v_{2},$ $\dots, v_{s})$ and $C_{2}(w_{1}, w_{2},\dots w_{l})$ be three cycles of type $A_{1}$. This gives $r \not= s \not= l$. Without loss of generality we assume $r < s$ and $r \not= l$. Now there is a shortest path, say $v_{1}\rightarrow u_{1}\rightarrow w_{1}$ of type $A_{2}$ or $A_{3}$ between $v_{1}$ and $w_{1}$ via $u_{1}$. This gives a path $v_{i}\rightarrow u_{i}\rightarrow w_{i}$ of type $A_{2}$ or $A_{3}$ between $v_{i}$ and $w_{i}$ via $u_{i}$. Hence by the assumption $r < s$, there does not exist any path of length two between $v_{r+j}$ and $w_{r}$ via $u_{r}$ for some $j(>0)$. This gives the link of $u_{r}$ is different from link of $u_{r-1}$. Which is a contradiction as the map $M$ of type $\{3^{3}, 4^{2}\}$. Hence $r = s = t$, that is, $length(C) = length(C_{1}) = length(C_{2})$. Similar argument we use for the cycle of type $A_{j}$ for $j \in \{2, 3\}$.
\hfill$\Box$\\

The cycles $C$ and $C_{t}$ for $t \in \{1, 2\}$ are homologous. Let $C_{1}, C_{2}, \dots, C_{m}$ be the possible list of homologous cycles of $C$. We claim that all of them have same length, that is,

\begin{lemma} \label{lem1.3} Let $C_{1}, C_{2}, \dots, C_{m}$ be the possible list of homologous cycles of type $A_{t}$ for $t \in \{1, 2, 3\}$. Then $length(C_{i})=length(C_{j})$ $\forall$ $i,j\in \{1, 2, \dots, m\}$.
\end{lemma}
\noindent{\sc Proof of Lemma}\ref{lem1.3} We consider a cylinder $S_{C_{i}}$ of $C_{i}$ of type $A_{1}$. Then we get a sub sequence $\{C_{i}, C_{i+1},\dots,C_{j}\}$ of $\{C_{1}, C_{2},\dots,C_{m}\}$ such that if we take any consecutive three cycles then they form a cylinder. Therefore, by the lemma \ref{lem1.2}, $length(C_{i}) =  length(C_{i+1}) = \dots =length(C_{j})$. This is true for any sub sequence of $\{C_{1}, C_{2},\dots, C_{m}\}$. Therefore, $length(C_{i}) = length(C_{j})$ $\forall$ $i,j\in \{1,2,...,m\}$. Similar argument we use for the cycle of type $A_{j}$ for $j \in \{2, 3\}$.
\hfill$\Box$\\

We define a $(r, s, k)$-representation of a map $M$ of type $\{3^{3}, 4^{2}\}$ on the torus. We consider a vertex $v$. By  definition \ref{def1.1}, \ref{def1.2} and \ref{def1.3}, we have three cycles at $v$, say $L_{1}(v), L_{2}(v), L_{3}(v)$, of types $A_{1},$ $A_{2}$ and $A_{3}$ at $v$ respectively. We cut $M$ along the cycle $L_{1}$, see Figure 1, $L_{1} = C(v_{1},v_{2},\dots,v_{r})$. Then starting at $v$ make another cut along $L_{3}$ until reaching $L_{1}$ again for the first time where the starting adjacent face to the horizontal base cycle $L_{1}$ is $4$-gon. In Figure 1, $v = v_{1}$. Let it crosses, say $s$, number of cycles homologous to $L_{1}$. In this section, the number $s$ of cycles which are homologous to $L_{1}$ is equal to the length of the path. In Figure 1, we took second cut along the path $v_{1}\rightarrow w_{1}\rightarrow x_{1}\rightarrow z_{1}\rightarrow v_{k+1}$ which is part of $L_{3}$ and $s = 4$. Let $length(L_{1}) = r$. We start from $v$ and denote the sequence of vertices along $L_{1}$ by $v = v_{1}, v_{2},\dots,v_{r}$. Hence we get a representation of $M$ and denoted it by {\em $(r, s)$-representation}. In the representation, we get identification of vertical sides in the natural manner but the identification of the horizontal sides needs some shifting in $(r, s)$-representation so that a vertex in the lower(base) side is identified with a vertex in the upper side. Let $v_{k+1}$ be the starting vertex of the upper horizontal cycle. Then $k$ denote the length of $v_{k+1}$ from $v_{1}$ in lower horizontal cycle $L_{1}$. Therefore we represent the $(r, s)$-representation by {\em $(r, s, k)$-representation}. In Figure 1, the vertex $v_{k+1}$ is the starting vertex of the upper horizontal cycle $C(v_{k+1},v_{k+2},\dots,v_{k})$ and $k = length(v_{1}\rightarrow v_{2} \rightarrow \dots \rightarrow v_{_{k+1}})$. By this construction we can show that every map of type $\{3^{3}, 4^{2}\}$ has a $(r, s, k)$-representation, that is, $T(r, s, k)$ . Therefore we claim it as a theorem,

\vspace{2 cm}
\begin{picture}(0,0)(-10,50)
\setlength{\unitlength}{2.6mm}
\drawpolygon(5,5)(20,5)(20,25)(5,25)
\drawpolygon(25,5)(40,5)(40,25)(25,25)


\drawline[AHnb=0](10,5)(10,25)
\drawline[AHnb=0](15,5)(15,25)
\drawline[AHnb=0](20,5)(20,25)
\drawline[AHnb=0](25,5)(25,25)
\drawline[AHnb=0](30,5)(30,25)
\drawline[AHnb=0](35,5)(35,25)

\drawline[AHnb=0](5,10)(20,10)
\drawline[AHnb=0](25,10)(40,10)
\drawline[AHnb=0](5,15)(20,15)
\drawline[AHnb=0](25,10)(40,10)
\drawline[AHnb=0](5,20)(20,20)
\drawline[AHnb=0](5,10)(20,10)

\drawline[AHnb=0](5,20)(10,25)
\drawline[AHnb=0](10,20)(15,25)
\drawline[AHnb=0](15,20)(20,25)
\drawline[AHnb=0](25,20)(30,25)
\drawline[AHnb=0](30,20)(35,25)
\drawline[AHnb=0](35,20)(40,25)
\drawline[AHnb=0](25,20)(40,20)

\drawline[AHnb=0](5,10)(10,15)
\drawline[AHnb=0](10,10)(15,15)
\drawline[AHnb=0](15,10)(20,15)
\drawline[AHnb=0](25,10)(30,15)
\drawline[AHnb=0](30,10)(35,15)
\drawline[AHnb=0](35,10)(40,15)
\drawline[AHnb=0](25,15)(40,15)

\put(21,10){\dots}
\put(21,5){\dots}
\put(21,15){\dots}
\put(21,20){\dots}
\put(21,25){\dots}

\put(5,4){$\scriptsize{v_{1}}$}
\put(10,4){$\scriptsize{v_{2}}$}
\put(15,4){$\scriptsize{v_{3}}$}
\put(20,4){$\scriptsize{v_{4}}$}
\put(25,4){$\scriptsize{v_{r-2}}$}
\put(30,4){$\scriptsize{v_{r-1}}$}
\put(35,4){$\scriptsize{v_{r}}$}
\put(40,4){$\scriptsize{v_{1}}$}

\put(5,9){$\scriptsize{w_{1}}$}
\put(10,9){$\scriptsize{w_{2}}$}
\put(15,9){$\scriptsize{w_{3}}$}
\put(20,9){$\scriptsize{w_{4}}$}
\put(25,9){$\scriptsize{w_{r-2}}$}
\put(30,9){$\scriptsize{w_{r-1}}$}
\put(35,9){$\scriptsize{w_{r}}$}
\put(40,9){$\scriptsize{w_{1}}$}

\put(5,14){$\scriptsize{x_{1}}$}
\put(10,14){$\scriptsize{x_{2}}$}
\put(15,14){$\scriptsize{x_{3}}$}
\put(20,14){$\scriptsize{x_{4}}$}
\put(25,14){$\scriptsize{x_{r-2}}$}
\put(30,14){$\scriptsize{x_{r-1}}$}
\put(35,14){$\scriptsize{x_{r}}$}
\put(40,14){$\scriptsize{x_{1}}$}

\put(5,19){$\scriptsize{z_{1}}$}
\put(10,19){$\scriptsize{z_{2}}$}
\put(15,19){$\scriptsize{z_{3}}$}
\put(20,19){$\scriptsize{z_{4}}$}
\put(25,19){$\scriptsize{z_{r-2}}$}
\put(30,19){$\scriptsize{z_{r-1}}$}
\put(35,19){$\scriptsize{z_{r}}$}
\put(40,19){$\scriptsize{z_{1}}$}

\put(5,24){$\scriptsize{v_{k+1}}$}
\put(10,24){$\scriptsize{v_{k+2}}$}
\put(15,24){$\scriptsize{v_{k+3}}$}
\put(20,24){$\scriptsize{v_{k+4}}$}
\put(25,24){$\scriptsize{v_{k-2}}$}
\put(30,24){$\scriptsize{v_{k-1}}$}
\put(35,24){$\scriptsize{v_{k}}$}
\put(40,24){$\scriptsize{v_{k+1}}$}

\put(13,0){\scriptsize Figure 1 : $T(r,4,k)$}

\end{picture}

\vspace{5.5 cm}

\begin{theo} \label{thm1.2} Every map of type $\{3^{3}, 4^{2}\}$ on the torus has a $(r, s, k)$-representation.
\end{theo}

\noindent{\sc Proof of Theorem}\ref{thm1.2} Proof of the theorem \ref{thm1.2} follows from the following argument.
\hfill$\Box$\\

Now we show the main result of this section that every map of type $\{3^{3}, 4^{2}\}$ on the torus has a contractible Hamiltonian cycle. Therefore we claim that,\

\begin{theo} \label{thm1.3} Let $M$ be a map of type $\{3^{3}, 4^{2}\}$ on the torus has a $(r, s, k)$-representation. Then it has a contractible Hamiltonian cycle.
\end{theo}

\noindent{\sc Proof of Theorem}\ref{thm1.3} By theorem \ref{thm1.2}, the map $M$ on the torus has a representation $T(r, s, k)$. The $T(r, s, k)$ has a rectangular planar shape with identical horizontal cycles and vertical paths on the boundary. Let $C_{1}, C_{2},...,C_{s}$ be the list of horizontal cycles of type $A_{1}$ (by the theorem \ref{thm1.1}) of length, say $r$ (by the lemma \ref{lem1.3}). Now, at each vertex we have only one cycle of type $A_{1}$. Therefore the cycles $C_{1}, C_{2},...,C_{s}$ will cover all the vertices of $M$, that is, $\bigcup_{i=1}^{s}V(C_{i}) = V(M)$. We define cycles as follows, $C_{1}(u_{1_{1}}, u_{1_{2}},\dots, u_{1_{r}})$, $C_{2}(u_{2_{1}}, u_{2_{2}},\dots, u_{2_{r}})$, $C_{3}(u_{3_{1}},u_{3_{2}},\dots, u_{3_{r}})$, \dots, $C_{s}(u_{s_{1}}, u_{s_{2}},\dots u_{s_{r}})$ in order, that is, if we consider any consecutive three cycles $C_{k},$ $C_{k+1},$ $C_{k+2}$ then $\partial S_{C_{k+1}} = \{C_{k}, C_{k+2}\}$.
In Figure 1, we have four cycles $C_{1}= C(u_{1}, u_{2},\dots, u_{r})$, $C_{2}= C(w_{1}, w_{2},\dots, w_{r})$, $C_{3}= C(x_{1}, x_{2},\dots, x_{r})$, and $C_{4}= C(z_{1}, z_{2},\dots, z_{r})$. We consider a cylinder which is bounded by $C_{i-1}$ and $C_{i}$ where $i$ is even integer and contains all quadrangles. For an even integer $i$, we define a new cycle $C_{i-1,i} := C(u_{i-1_{1}},$ $u_{i-1_{2}},$ $\dots, u_{i-1_{r}},$ $u_{i_{r}},$ $u_{i_{r-1}},\dots,$ $u_{i_{1}})$. It is contractible as it bounds a $2$-disk which contains only quadrangles. In Figure 1, for $i=2$ the cycle is $C(v_{1},$ $v_{2},$ $\dots,$ $v_{r},$ $w_{r},$ $w_{r-1},$ $\dots,$ $w_{1})$ which is defined by $C_{1}$ and $C_{2}$. Therefore, for each pair of cycles $\{C_{i-1}, C_{i}\}$ where $i = 2, 4,\dots,s$, we have a contractible cycle $C_{i-1, i}$. We consider two cycles $C_{i-1, i}$, $C_{i+1, i+2}$ and adjacent two $3$-gons $\triangle(u_{i_{1}}, u_{i_{2}}, u_{i+1_{2}})$ and $\triangle(u_{i+1_{1}}, u_{i+1_{2}}, u_{i_{1}})$. For an even integer $i$, we define a cycle $C_{i-1, i}\cup C_{i+1, i+2} := C($ $u_{i-1_{1}},$, $u_{i-1_{2}},$ \dots $u_{i-1_{r}},$ $u_{i_{r}},$ $u_{i_{r-1}},$ \dots $u_{i_{2}},$ $u_{i+1_{2}},$, $u_{i+1_{3}},$ \dots, $u_{i+1_{r}},$ $u_{i+2_{r}},$ $u_{i+2_{r-1}},$ \dots, $u_{i+2_{1}},$ $u_{i+1_{1}},$ $u_{i_{1}})$ which is concatenation of $C_{i-1, i}$ and $C_{i+1,i+2}$. Therefore we define cycle $C \colon = C_{1, 2}\cup C_{3, 4} \cup \dots \cup C_{s-1, s}$. The cycle $C$ is contractible as each $C_{i-1, i}$ is contractible. Also the cycle $C$ is Hamiltonian as it contains all the vertices of $C_{i}$ for $i \in \{1,\dots,s\}$. Therefore the map $M$ is Hamiltonian and contains contractible Hamiltonian cycle.
\hfill$\Box$

\section{ Map of type $\{3^{2}, 4, 3, 4\}$}

We consider a map $M$ of type $\{3^{2}, 4, 3, 4\}$. We define a path in $M$ as follows $\colon$

\begin{defn}\label{df2.1} Let $P(u_{1}\rightarrow\dots \rightarrow u_{i-1}\rightarrow u_{i} \rightarrow u_{i+1} \rightarrow \dots \rightarrow u_{r})$ be a path in $M$. The path $P$ is said to be well defined at a vertex $u_{i}$ if we denote $u_{i+1}$, $u_{i-1}$ be the successor and predecessor of $u_{i}$ respectively in $P$ or in the extended path of $P$ then,\
     \begin{enumerate}
        \item the link of $u_{i}$, $lk(u_{i})=C(\textbf{a},$ $u_{i+1}, b,$ $c, \textbf{d}, u_{i-1}, e)$ gives $lk(u_{i-1})=C(\textbf{f}, g,$ $e, u_{i}, \textbf{c},$ $d, u_{i-2})$ and $lk(u_{i+1})=C(\textbf{e},$ $a, k, u_{i+2}, \textbf{l},$ $b, u_{i})$ and

        \item the link of $u_{i}$, $lk(u_{i})=C(\textbf{e},$ $h, k, u_{i+1}, \textbf{l},$ $b, u_{i-1})$ gives $lk(u_{i-1})=C(\textbf{h}, u_{i},$ $b, c, \textbf{d},$ $u_{i-2}, e)$ and $lk(u_{i+1})=C(\textbf{s},$ $u_{i+2}, t, l, $\textbf{b},$ u_{i}, k)$.

     \end{enumerate}

We denote this path by $B_{1}$.
\end{defn}

We consider a maximal path $P$ of type $B_{1}$. The path $P$ is cycle by the similar argument of theorem \ref{thm1.1}. We consider two cycles $C_{1}$ and $C_{2}$ of type $B_{1}$. We claim the next lemma ,\

\begin{lemma}\label{lem2.1} Let $C_{1}$ and $C_{2}$ be two cycles of type $B_{1}$. Then $C_{1} = C_{2}$ if $E(C_{1}) \cap E(C_{2})\not= \emptyset$.
\end{lemma}

\noindent{\sc Proof of Lemma}\ref{lem2.1} Let $C_{1}(u_{1_{1}}, u_{1_{2}},\dots, u_{1_{r}})$ and $C_{2}(u_{2_{1}}, u_{2_{2}},\dots, u_{2_{s}})$ be two cycles. Let $E(C_{1}) \cap E(C_{2}) \not= \emptyset$ then there is a $e \in E(C_{1} \cap C_{2})$ where $e = yx$. The cycles $C_{1}, C_{2}$ both are well defined at {\em y, x}. We consider the link of $x$. If $lk(x) = C(\textbf{a}, b, c, w, \textbf{d}, e, y)$ then $ w \in V(C_{i})$ for $i \in \{1, 2\}$ by the definition \ref{df2.1}. This gives, the path $y\rightarrow x\rightarrow w$ is  part of both $C_{1}$ and $C_{2}$, that is, $y = u_{1_{t_{1}-1}} = u_{2_{t_{2}-1}}$, $x = u_{1_{t_{1}}} = u_{2_{t_{2}}}$ and $w = u_{1_{t_{1}+1}} = u_{2_{t_{2}+1}}$ for some $t_{1}\in \{1,\dots,r\}$ and $t_{2} \in \{1,\dots,s\}$. Again, we consider the edge $xw$ and continue. We stop after $r$ steps. Let $t_{2}-t_{1} = m$ for some $m$. Hence we get $u_{1_{1}} = u_{2_{m+1}}, u_{1_{2}} = u_{2_{m+2}},\dots, u_{1_{r}} = u_{2_{m+r}}$ and $u_{1_{1}} = u_{2_{m+r+1}}$. This imply that $m+1 = m+r+1$ and $r = s$ as $u_{1_{m+r+1}} = u_{2_{m+1}}$ and $C_{2}$ is cycle. This gives the cycles are identical. If $lk(x) = C(\textbf{a}, b, c, w, \textbf{d}, y, e)$ then $ b \in V(C_{i})$ for $i \in \{1, 2\}$. By the similar argument, the cycles $C_{1}$ and $C_{2}$ are identical. Therefore, $E(C_{1}) \cap E(C_{2})\not= \emptyset$ gives $C_{1} = C_{2}$.
\hfill$\Box$\\

We choose a vertex $u$ in $C$ of type $B_{1}$. We collect all adjacent faces of $u$ in a set, say $S$. We do it for all of $V(C)$. We consider geometric carrier of $S$ which is a cylinder as the boundaries of the cylinder are either two disjoint cycles or an identical cycle. We denote it by $S_{C}$. Let $\partial S_{C} = \{C_{1}, C_{2}\}$. Now we claim the next lemma,\

\begin{lemma}\label{lem3.2} Let $C$ be a cycle of type $B_{1}$ and $\partial S_{C} = \{C_{1}, C_{2}\}$. Then $C_{1} = C_{2}$ if $C_{1} \cap C_{2}\not= \emptyset$.
\end{lemma}

\noindent{\sc Proof of Lemma}\ref{lem3.2} Suppose $C_{1} \cap C_{2}$ is non empty then $C_{1} \cap C_{2}$ is a set of vertices or contains an edge. Now we have the following cases,\
\begin{enumerate}
\item Let $C_{1} \cap C_{2}$ be a set of vertices and does not contain edge. Let $u \in V(C_{1} \cap C_{2})$. We consider adjacent edges of $u$ where $u \in V(C_{i})$ for $i \in \{1, 2\}$. Now by definition \ref{df2.1}, the number of adjacent edges on one side of the cycle $C_{i}$ is two and on the other side one at each vertex. Therefore, the $degree(u)\geq3$ or $4$ in the cylinder for $u \in V(C_{i})$. If we assume $degree(u)\geq3$ for $u \in V(C_{1})$ then, for $u \in V(C_{2})$ $degree(u)\geq4$ in cylinder. This gives $degree(u)\geq7(=3+4)$. Hence the degree of $u$ is at least seven which is a contradiction as the degree of $u$ is five. Therefore $C_{1} \cap C_{2}$ contains an edge. Now, by the lemma \ref{lem2.1}, $C_{1} = C_{2}$.

\item If $C_{1} \cap C_{2}$ contains an edge then by lemma \ref{lem2.1}, $C_{1} = C_{2}$.
\hfill$\Box$\\
\end{enumerate}

\vspace{1.5 cm}
\begin{picture}(0,0)(-10,50)
\setlength{\unitlength}{2.6 mm}
\drawpolygon(5,5)(20,5)(20,25)(5,25)
\drawpolygon(25,5)(40,5)(40,25)(25,25)


\drawline[AHnb=0](10,5)(10,25)
\drawline[AHnb=0](15,5)(15,25)
\drawline[AHnb=0](20,5)(20,25)
\drawline[AHnb=0](25,5)(25,25)
\drawline[AHnb=0](30,5)(30,25)
\drawline[AHnb=0](35,5)(35,25)

\drawline[AHnb=0](5,10)(20,10)
\drawline[AHnb=0](25,10)(40,10)
\drawline[AHnb=0](5,15)(20,15)
\drawline[AHnb=0](25,10)(40,10)
\drawline[AHnb=0](5,20)(20,20)
\drawline[AHnb=0](5,10)(20,10)

\drawline[AHnb=0](5,20)(10,25)
\drawline[AHnb=0](10,20)(15,15)
\drawline[AHnb=0](15,20)(20,25)
\drawline[AHnb=0](25,20)(30,15)
\drawline[AHnb=0](30,20)(35,25)
\drawline[AHnb=0](35,20)(40,15)
\drawline[AHnb=0](25,20)(40,20)

\drawline[AHnb=0](5,10)(10,15)
\drawline[AHnb=0](10,10)(15,5)
\drawline[AHnb=0](15,10)(20,15)
\drawline[AHnb=0](25,10)(30,5)
\drawline[AHnb=0](30,10)(35,15)
\drawline[AHnb=0](35,10)(40,5)
\drawline[AHnb=0](25,15)(40,15)

\put(21,10){\dots}
\put(21,5){\dots}
\put(21,15){\dots}
\put(21,20){\dots}
\put(21,25){\dots}

\put(5,4){$\scriptsize{v_{1}}$}
\put(10,4){$\scriptsize{v_{2}}$}
\put(15,4){$\scriptsize{v_{3}}$}
\put(20,4){$\scriptsize{v_{4}}$}
\put(25,4){$\scriptsize{v_{r-2}}$}
\put(30,4){$\scriptsize{v_{r-1}}$}
\put(35,4){$\scriptsize{v_{r}}$}
\put(40,4){$\scriptsize{v_{1}}$}

\put(5,9){$\scriptsize{w_{1}}$}
\put(10,10.3){$\scriptsize{w_{2}}$}
\put(15,9){$\scriptsize{w_{3}}$}
\put(20,9){$\scriptsize{w_{4}}$}
\put(25,10.3){$\scriptsize{w_{r-2}}$}
\put(30,9){$\scriptsize{w_{r-1}}$}
\put(35,10.3){$\scriptsize{w_{r}}$}
\put(40,9){$\scriptsize{w_{1}}$}

\put(5,14){$\scriptsize{x_{1}}$}
\put(10,14){$\scriptsize{x_{2}}$}
\put(15,14){$\scriptsize{x_{3}}$}
\put(20,14){$\scriptsize{x_{4}}$}
\put(25,14){$\scriptsize{x_{r-2}}$}
\put(30,14){$\scriptsize{x_{r-1}}$}
\put(35,14){$\scriptsize{x_{r}}$}
\put(40,14){$\scriptsize{x_{1}}$}

\put(5,19){$\scriptsize{z_{1}}$}
\put(10,20.3){$\scriptsize{z_{2}}$}
\put(15,19){$\scriptsize{z_{3}}$}
\put(20,19){$\scriptsize{z_{4}}$}
\put(25,20.3){$\scriptsize{z_{r-2}}$}
\put(30,19){$\scriptsize{z_{r-1}}$}
\put(35,20.3){$\scriptsize{z_{r}}$}
\put(40,19){$\scriptsize{z_{1}}$}

\put(5,24){$\scriptsize{v_{2k+1}}$}
\put(10,24){$\scriptsize{v_{2k+2}}$}
\put(15,24){$\scriptsize{v_{2k+3}}$}
\put(20,24){$\scriptsize{v_{2k+4}}$}
\put(25,24){$\scriptsize{v_{2k-2}}$}
\put(30,24){$\scriptsize{v_{2k-1}}$}
\put(35,24){$\scriptsize{v_{2k}}$}
\put(40,24){$\scriptsize{v_{2k+1}}$}

\put(14,0) {\scriptsize Figure 2: $T(r, 4, 2k)$}

\end{picture}

\vspace{5.5 cm}

Therefore the boundaries of a cylinder $S_{C}$ are either two disjoint cycles or an identical cycle. Hence $S_{C}$ is a cylinder with two boundaries $C_{1}$ and $C_{2}$ of type $B_{1}$. We consider a cylinder $S_{C_{i}}$ which is bounded by two cycles, say $C_{l},C_{k}$. Then they are of same length by the similar argument of lemma \ref{lem1.2}. If we consider any two homologous cycles $C_{1}, C_{2}$ of type $B_{1}$ then $length(C_{1}) = length(C_{2})$ by the similar argument of lemma \ref{lem1.3}. By the definition \ref{df2.1} we have two cycles of type $B_{1}$ at a vertex $v$, say $L_{1}(v), L_{2}(v)$. We define $T(r, s, k)$ of $M$ by the construction of $(r,s,k)$-representation done in Section 2. Here, we take the second cut along the cycle $L_{2}(v)$ where the starting adjacent face to the base horizontal cycle is $4$-gon. Hence we get $T(r, s, k)$ of the map $M$ of type $\{3^{2}, 4, 3, 4\}$. This $T(r, s, k)$ exists for any map of type $\{3^{2}, 4, 3, 4\}$ on the torus. See Figure 2, it is $T(r, 4, 2k)$. Here $L_{1}=C(v_{1}, v_{2}, \dots, v_{r})$ and the path $v_{1}\rightarrow w_{1}\rightarrow x_{1}\rightarrow z_{1}\rightarrow v_{2k+1}$ which is part of the cycle $L_{2}$. Along this path we took the second cut.  We consider $C_{1}, C_{2},\dots,C_{s}$ be the sequence of all possible homologous horizontal cycles of length, say $r$ ( by the lemma \ref{lem2.1}) of type $B_{1}$ in $T(r, s, k)$. We use similar argument like done in theorem \ref{thm1.3}. Hence we get a cycle $C$ which is contractible and Hamiltonian. Therefore the map of type $\{3^{2}, 4, 3, 4\}$ on the torus is Hamiltonian.

\begin{theo} \label{thm2.1} Let $M$ be a map of type $\{3^{2}, 4, 3, 4\}$ on the torus. Then it is Hamiltonian.
\end{theo}

\noindent{\sc Proof of Theorem}\ref{thm2.1} The proof of the theorem \ref{thm2.1} follows from the result of this section.
\hfill$\Box$

\section{Map of type $\{ 6, 3, 6, 3\}$}

We consider a map $M$ of type $\{ 6, 3, 6, 3\}$. We define a path as follows $\colon$

\begin{defn}\label{df3.1} Let $P(u_{1}\rightarrow u_{r})$ be a path in $M$. We denote $A_{v}$ be a set of edges which are adjacent at a vertex $v$ in $M$. The path $P$ is said to be well defined at an inner vertex $u_{i}$ if,\
  \begin{enumerate}
  \item two edges of $M$ emerge from $u_{i}$ which are not in $E(P)\cap A_{u_{i}}$,
  \item one of them on one side of $u_{i-1}\rightarrow u_{i}\rightarrow u_{i+1}$and remaining one on the other side.
  \item The boundary vertex of $P$ also well defined in the extended path of $P$, that is, there is a path which contains the path $P$ and satisfy $(1)$ and $(2)$ and both the boundary vertices of $P$ are inner vertices.
  \end{enumerate}
We denote this path by $X_{1}$.
\end{defn}

\vspace{1 mm}
\begin{picture}(0,0)(-10,50)
\setlength{\unitlength}{2.4mm}

\drawline[AHnb=0](0,0)(40,0)
\drawline[AHnb=0](5,10)(45,10)
\drawline[AHnb=0](10,20)(50,20)

\drawline[AHnb=0](0,0)(10,20)
\drawline[AHnb=0](10,0)(20,20)
\drawline[AHnb=0](20,0)(30,20)
\drawline[AHnb=0](30,0)(40,20)
\drawline[AHnb=0](40,0)(50,20)

\drawline[AHnb=0](2.5,5)(5,0)
\drawline[AHnb=0](7.5,15)(15,0)
\drawline[AHnb=0](15,20)(25,0)
\drawline[AHnb=0](25,20)(35,0)
\drawline[AHnb=0](35,20)(42.5,5)
\drawline[AHnb=0](45,20)(47.5,15)

\put(0,-1){$\scriptsize{v_{1}}$}
\put(5,-1){$\scriptsize{v_{2}}$}
\put(10,-1){$\scriptsize{v_{3}}$}
\put(15,-1){$\scriptsize{v_{4}}$}
\put(20,-1){$\scriptsize{v_{5}}$}
\put(25,-1){$\scriptsize{v_{6}}$}
\put(30,-1){$\scriptsize{v_{7}}$}
\put(35,-1){$\scriptsize{v_{8}}$}
\put(40,-1){$\scriptsize{v_{1}}$}

\put(3,5){$\scriptsize{v_{9}}$}
\put(13,5){$\scriptsize{v_{10}}$}
\put(23,5){$\scriptsize{v_{11}}$}
\put(33,5){$\scriptsize{v_{12}}$}
\put(43,5){$\scriptsize{v_{9}}$}

\put(5,9){$\scriptsize{v_{13}}$}
\put(10.3,9){$\scriptsize{v_{14}}$}
\put(15,9){$\scriptsize{v_{15}}$}
\put(20.3,9){$\scriptsize{v_{16}}$}
\put(25,9){$\scriptsize{v_{17}}$}
\put(30.3,9){$\scriptsize{v_{18}}$}
\put(35,9){$\scriptsize{v_{19}}$}
\put(40.3,9){$\scriptsize{v_{20}}$}
\put(45,9){$\scriptsize{v_{13}}$}

\put(8,15){$\scriptsize{v_{21}}$}
\put(18,15){$\scriptsize{v_{22}}$}
\put(28,15){$\scriptsize{v_{23}}$}
\put(38,15){$\scriptsize{v_{24}}$}
\put(48,15){$\scriptsize{v_{21}}$}

\put(10,20.5){$\scriptsize{v_{7}}$}
\put(15,20.5){$\scriptsize{v_{8}}$}
\put(20,20.5){$\scriptsize{v_{1}}$}
\put(25,20.5){$\scriptsize{v_{2}}$}
\put(30,20.5){$\scriptsize{v_{3}}$}
\put(35,20.5){$\scriptsize{v_{4}}$}
\put(40,20.5){$\scriptsize{v_{5}}$}
\put(45,20.5){$\scriptsize{v_{6}}$}
\put(50,20.5){$\scriptsize{v_{7}}$}

\put(15,-5) {\scriptsize Figure 3 : $T(8, 2, 6)$}

\end{picture}

\vspace{7 cm}

We consider a maximal path of type $X_{1}$ which is cycle by the similar argument of theorem \ref{thm1.1}. By definition \ref{df3.1}, we can define two cycles of type $X_{1}$ at each vertex. We consider two cycles of type $X_{1}$ at a vertex $v$, say $L_{1}(v), L_{2}(v)$. We define $(r, s, k)$-representation of a map $M$ by the similar argument given in Section 2. Hence we get $T(r, s, k)$ of $M$. Figure 3 is a $T(8, 2, 6)$ of a map of type $\{ 6, 3, 6, 3\}$ on 24 vertices. In Figure 3, $L_{1} = C(v_{1},$ $v_{2},$ $v_{3},$ $\dots,$ $v_{8})$ and the path $v_{1}\rightarrow v_{9}\rightarrow v_{13}\rightarrow v_{21}\rightarrow v_{7}$ which is part of $L_{2}$.  In the process we take second cut along the cycle $L_{2}$ where the starting adjacent face to the base horizontal cycle $L_{1}$ is $3$-gon. Now we claim the following theorem,\


\begin{theo} \label{thm3.2} Let $M$ be a map of type $\{ 6, 3, 6, 3\}$ on the torus which has a $T(r, s, k)$. Then it contains a contractible Hamiltonian cycle.
\end{theo}

\noindent{\sc Proof of Theorem}\ref{thm3.2} We denote $C_{1}$ be horizontal base cycle of $T(r, s, k)$. We consider all horizontal cycles which are homologous to $C_{1}$, say $C_{2}, C_{3},\dots,C_{t}$ of length, say $r$ (by the similar argument of lemma \ref{lem1.2}). Let $C_{1}(u_{1,1},$ $u_{1,2},$ $\dots,$ $u_{1,r})$, $C_{2}(u_{2,1},$ $u_{2,2},$ $\dots,$ $u_{2,r})$,$ $\dots,  $C_{t}(u_{t,1},$ $u_{t,2},$ $\dots,$ $u_{t,r})$ be the cycles. In Figure 3, $C_{1}= C(v_{1},$ $\dots,$ $v_{8})$ and $C_{2}= C(v_{13},$ $\dots,$ $v_{20})$. We denote a vertex $u_{(i_{k-1},i_{k}),(i+1_{k-2},i+1_{k-1})}$ which does not belongs to any one of the cycles $C_{1}, C_{2},\dots, C_{t}$. Let $u_{(i_{k-1},i_{k}),(i+1_{k-2},i+1_{k-1})}$ adjacent to the vertices $u_{i_{k-1}}$, $u_{i_{k}}$ of $C_{i}$ and $u_{i+1_{k-2}}$, $u_{i+1_{k-1}}$ of $C_{i+1}$. Also there are two $3$-gons $\triangle(u_{i_{k-1}}$ $\rightarrow$ $ u_{i_{k}}$ $\rightarrow$ $u_{(i_{k-1},i_{k}),(i+1_{k-2},i+1_{k-1})}$ $\rightarrow$ $ u_{i_{k-1}})$ and $\triangle(u_{i+1_{k-2}}$ $\rightarrow$ $ u_{i+1_{k-1}}$ $\rightarrow$ $ u_{(i_{k-1},i_{k}),(i+1_{k-2},i+1_{k-1})}$ $\rightarrow$ $ u_{i+1_{k-2}})$ which are adjacent to the vertex $u_{(i_{k-1},i_{k}),(i+1_{k-2},i+1_{k-1})}$. In Figure 3, $t = 2$ and $v_{(1_{5},1_{6}),(2_{16},2_{17})} = v_{11}$. When $t = 1,$ we use $u_{(i_{k-1},i_{k})}$ in place of $u_{(i_{k-1},i_{k}),(i+1_{k-2},i+1_{k-1})}$ as upper and base horizontal cycles are same in $T(r,s,k)$. Now we follow one of the following cases$\colon$

\begin{enumerate}
\item Let $t = 1$. We consider $C_{1}(u_{1,1}, u_{1,2}, \dots, u_{1,r})$. Also we consider adjacent $3$-gons to $C_{1}$. We concatenate the cycle $C_{1}$ and $3$-gons which are adjacent to the cycle $C_{1}$ and contains a new vertex. Hence we get a cycle $C(u_{1,1}, u_{(1_{1},1_{2})},$ $u_{1,2},$ $u_{1,3},$ $u_{(1_{3},1_{4})},$ $u_{1,4},$ $\dots,$ $u_{1,r-1},$ $u_{(1_{r-1},1_{r})}, u_{1,r})$. This cycle contains all the vertices of $M$.

\item Let $t = 2$. We denote a set $S := \emptyset$. We consider all adjacent $6$-gons to the cycle $C_{1}$ in $S$ except last one. In Figure 3, last $6$-gon is $F_{6}(v_{1}, v_{9}, v_{20}, v_{19}, v_{12}, v_{8})$. Also we consider all adjacent $3$-gons to $C_{1}$ in $S$. The geometric carrier $|S|$ is a $2$-disk. If there is a vertex which is not on the boundary of $2$-disk $|S|$ then there is a $3$-gon adjacent along an edge to $2$-disk and contains the vertex. We put this $3$-gon in $S$. Hence we get $S$ where the geometric carrier $|S|$ is $2$-disk. The boundary $\partial|S|$ is a cycle which contains all the vertices of $M$. Therefore, we get a cycle which contains all the vertices of $M$. In Figure 3, the cycle $\partial|S|$ is $C(v_{1},$ $v_{9},$ $v_{20},$ $v_{13},$ $v_{21},$ $v_{14},$ $v_{10},$ $v_{15},$ $v_{22},$ $v_{16},$ $v_{11},$ $v_{17},$ $v_{23},$ $v_{18},$ $v_{19},$ $v_{12},$ $v_{8},$ $v_{7},$ $v_{6},$ $v_{5},$ $v_{24},$ $v_{4},$ $v_{3},$ $v_{2})$.

\item Let $t > 2$ and $t$ is an even integer. We consider the cycles $C_{i}$ and $C_{i+1}$ where $i$ is an odd integer. Now we use similar idea of $(2)$. Hence we get a $2$-disk, say $D_{i}$. We do it for all odd integer $i$ where $i \in \{1,\dots,t\}$. Hence we get $2$-disks $D_{1}$, $D_{3}$, \dots, $D_{t-1}$. Let $D_{i}$ and $D_{i+2}$ be two $2$-disks where $i$ is an odd integer. Let $F_{1}, F_{2}$ be two faces. We denote $|\{F_{1},F_{2}\}|$ is geometric carrier of the set $\{F_{1}, F_{_{2}}\}$. Now there are two faces $6$-gon, say $F_{i,6}$ and $3$-gon, say $F_{i,3}$ between $D_{i}$ and $D_{i+2}$ such that $E(|\{F_{i, 3},F_{i, 6}\}|)\cap E(D_{t})$ is an edge for $t \in \{i, i+2\}$. The geometric carrier of $\{ F_{i, 3}, F_{i, 6}, D_{i}, D_{i+2}\}$ is a $2$-disk, say $D_{i,i+2}$. Also $V(D_{i,i+2}) = V(D_{i}) \cup V(D_{i+2})$. We do it for any two consecutive $2$-disks $D_{i}$, $D_{i+2}$ for $i \in \{1, 3, 5, \dots, t-1\}$. Hence we get a $2$-disk, say $D:=\cup_{i}D_{i,i+2}$. Which is contractible and contains all the vertices of $M$.

\item Let $t > 2$  and $t$ is an odd integer. We follow $(3)$ for $C_{1}, C_{2},\dots,C_{t-1}$. We get a $2$-disk, say $D$. Next we follow $(1)$ with the cycle $C_{t}$ in place of $C_{1}$. Hence we get a cycle $C$ which contains all the vertices of $C_{t}$ and the vertices other than $V(D)$ which are belongs to $3$-gons and adjacent to $C_{t}$. Now there are two faces $6$-gon, say $f_{6}$, and $3$-gon, say $f_{3}$, between $C$ and $D$ such that $E(|\{f_{3},f_{6}\}|)\cap E(D)$ and $E(|\{f_{3},f_{6}\}|)\cap E(C)$ are two edges. This gives a cycle which is concatenation of $\partial D$ and $C$ by $\partial|\{f_{3},f_{6}\}|$. Hence we get a cycle which contains all the vertices of $M$.
\end{enumerate}
Therefore, the map of type $\{ 6, 3, 6, 3\}$ is Hamiltonian
\hfill$\Box$

\section{Map of type $\{ 3^{4}, 6\}$}

Let $M$ be a map of type $\{ 3^{4}, 6\}$. We define paths as follows $\colon$

\begin{defn}\label{df4.1} Let $P(u_{1}\rightarrow u_{r})$ be a path in $M$. We denote three paths $Q_{1}(w_{i}\rightarrow w_{i+1}\rightarrow w_{i+2} \rightarrow w_{i+3})$ (see Figure 4),  $Q_{2}(v_{i}\rightarrow v_{i+1}\rightarrow v_{i+2} \rightarrow v_{i+3})$ (see Figure 5) and $Q_{3}(z_{i}\rightarrow z_{i+1}\rightarrow z_{i+2} \rightarrow z_{i+3})$ where $link(w_{i})=C(w_{i-1},$ $x_{2},w_{i+1},$ $\textbf{w}_{i+2},$ $\textbf{x}_{3},$ $\textbf{x}_{4},$ $ x_{5},$ $x_{6})$, $link(w_{i+1})=C(w_{i},$ $x_{2},$ $x_{7},$ $x_{9},$ $w_{i+2},$ $\textbf{x}_{3},$ $ \textbf{x}_{4},$ $\textbf{x}_{5})$, $link(w_{i+2})=C(w_{i+1},$ $x_{8},$ $x_{9},$ $w_{i+3},$ $x_{3},$ $\textbf{x}_{4},$ $ \textbf{x}_{5},$ $\textbf{w}_{i},)$, $link(w_{i+3})=C(w_{i+2},$ $x_{9},$ $w_{i+4},$ $\textbf{w}_{i+5},$ $\textbf{x}_{10},$ $\textbf{x}_{11},$ $ x_{12},x_{3})$, $link(v_{i})=C(v_{i-1},$ $y_{1},y_{2},$ $y_{3},$ $v_{i+1},$ $\textbf{y}_{4},$ $ \textbf{y}_{5},$ $\textbf{y}_{6})$, $link(v_{i+1})=C(v_{i},$ $y_{3},$ $y_{7},$ $v_{i+2},$ $y_{4},$ $\textbf{y}_{5},$ $ \textbf{y}_{6},$ $\textbf{v}_{i-1})$, $link(v_{i+2})=C(v_{i+1},$ $v_{7},$ $v_{i+3},$ $\textbf{v}_{i+4},$ $\textbf{y}_{11},$ $\textbf{y}_{12},$ $ y_{8},$ $w_{4},)$, $link(v_{i+3})=C(v_{i+2},$ $y_{7},$ $y_{9},$ $y_{10},$ $v_{i+4},$ $\textbf{y}_{11},$ $ \textbf{y}_{12},\textbf{y}_{8})$, $link(z_{i})=C(z_{i-1},$ $r_{3},r_{4},$ $z_{i+1},$ $r_{12},$ $ \textbf{r}_{13},$ $\textbf{r}_{14}, \textbf{z}_{i-2})$, $link(z_{i+1})=C(z_{i},$ $r_{4},$ $z_{i+2},$ $\textbf{z}_{i+3},$ $\textbf{r}_{9},$ $\textbf{r}_{10},$ $ \textbf{r}_{11},$ $r_{12})$, $link(z_{i+2})=C(z_{i+1},$ $r_{4},$ $r_{5},$ $r_{6},$ $z_{i+3},$ $\textbf{r}_{9},$ $ \textbf{r}_{10},$ $\textbf{r}_{11},)$ and $link(z_{i+3})=C(z_{i+2},$ $r_{6},$ $r_{7},$ $z_{i+4},$ $r_{9},$ $\textbf{r}_{10},$ $ \textbf{r}_{11},\textbf{z}_{i+1})$. A path $P$ is said to be well defined if we consider a path of length three, say $L(u_{t}\rightarrow u_{t+1}\rightarrow u_{t+2}\rightarrow u_{t+3})$ in P or in the extended path of $P$ then either $L \mapsto Q_{1}$ by $u_{j}\mapsto w_{j}$ or $L \mapsto Q_{2}$ by $u_{j}\mapsto v_{j}$ or $L \mapsto Q_{3}$ by $u_{j}\mapsto z_{j}$ for $j \in \{t,t+1,t+2,t+3\}$. In Figure 6, $P = path(u_{1}\rightarrow u_{2}\rightarrow \dots \rightarrow u_{10})$. We denote this path by $Y_{1}$.
\end{defn}

\vspace{1 mm}
\begin{picture}(0,0)(-30,20)
\setlength{\unitlength}{1.5 mm}

\drawpolygon(-15,0)(-10,-5)(-5,0)(-5,5)(-5,10)(-10,15)(-15,10)
\drawpolygon(-5,0)(5,0)(5,10)(-5,10)
\drawpolygon(5,0)(10,-5)(15,0)(15,5)(15,10)(10,15)(5,10)

\drawpolygon(-15,0)(-10,-5)(-10,-15)(-15,-10)
\drawpolygon(-5,0)(5,0)(5,-10)(-5,-10)
\drawpolygon(-5,0)(-10,-5)(-10,-15)(-5,-10)

\drawpolygon(5,0)(10,-5)(10,-15)(5,-10)
\drawpolygon(15,0)(10,-5)(10,-15)(15,-10)


\drawline[AHnb=0](-10,-5)(-15,-10)
\drawline[AHnb=0](-10,-5)(-5,-10)
\drawline[AHnb=0](-5,0)(5,-10)
\drawline[AHnb=0](5,0)(-5,10)

\drawline[AHnb=0](10,-5)(5,-10)
\drawline[AHnb=0](10,-5)(15,-10)


\put(-18,-10){$\scriptsize{y_{1}}$}
\put(-10,-17){$\scriptsize{y_{2}}$}
\put(-5,-12){$\scriptsize{y_{3}}$}
\put(3,-12){$\scriptsize{y_{7}}$}
\put(10,-17){$\scriptsize{y_{9}}$}
\put(15,-12){$\scriptsize{y_{10}}$}

\put(-20,0){$\scriptsize{v_{i-1}}$}
\put(-11,-3){$\scriptsize{v_{i}}$}
\put(-10,0){$\scriptsize{v_{i+1}}$}
\put(6,0){$\scriptsize{v_{i+2}}$}
\put(7.8,-3){$\scriptsize{v_{i+3}}$}
\put(15,0){$\scriptsize{v_{i+4}}$}

\put(-17,10){$\scriptsize{y_{6}}$}
\put(-11,12){$\scriptsize{y_{5}}$}
\put(-8,10){$\scriptsize{y_{4}}$}
\put(6,10){$\scriptsize{y_{8}}$}
\put(9,12){$\scriptsize{y_{12}}$}
\put(14,11){$\scriptsize{y_{11}}$}

\put(-19,-24){\scriptsize Figure 4 : $path(v_{i}\rightarrow v_{i+1}\rightarrow v_{i+2} \rightarrow v_{i+3})$}


\drawpolygon(25,0)(35,0)(35,10)(25,10)
\drawpolygon(35,0)(40,-5)(45,0)(45,5)(45,10)(40,15)(35,10)
\drawpolygon(45,0)(55,0)(55,10)(45,10)
\drawpolygon(55,0)(60,-5)(65,0)(65,5)(65,10)(60,15)(55,10)

\drawpolygon(25,0)(35,0)(35,-10)(25,-10)

\drawpolygon(35,0)(40,-5)(40,-15)(35,-10)
\drawpolygon(45,0)(55,0)(55,-10)(45,-10)
\drawpolygon(45,0)(40,-5)(40,-15)(45,-10)

\drawpolygon(55,0)(60,-5)(60,-15)(55,-10)
\drawpolygon(65,0)(60,-5)(60,-15)(65,-10)
%
\drawline[AHnb=0](25,0)(35,-10)
\drawline[AHnb=0](35,0)(25,10)

\drawline[AHnb=0](40,-5)(35,-10)
\drawline[AHnb=0](40,-5)(45,-10)
\drawline[AHnb=0](45,0)(55,-10)
\drawline[AHnb=0](55,0)(45,10)

\drawline[AHnb=0](60,-5)(55,-10)
\drawline[AHnb=0](60,-5)(65,-10)

\put(25,-12){$\scriptsize{x_{1}}$}
\put(34,-12){$\scriptsize{x_{2}}$}
\put(40,-17){$\scriptsize{x_{7}}$}
\put(45,-12){$\scriptsize{x_{8}}$}
\put(54,-12){$\scriptsize{x_{9}}$}
\put(60,-17){$\scriptsize{x_{13}}$}
\put(65,-12){$\scriptsize{x_{14}}$}

\put(20,0){$\scriptsize{w_{i-1}}$}
\put(36,0){$\scriptsize{w_{i}}$}
\put(37,-3){$\scriptsize{w_{i+1}}$}
\put(39,0){$\scriptsize{w_{i+2}}$}
\put(56,0){$\scriptsize{w_{i+3}}$}
\put(57,-3){$\scriptsize{w_{i+4}}$}
\put(66,0){$\scriptsize{w_{i+5}}$}

\put(24,11){$\scriptsize{x_{6}}$}
\put(33,11){$\scriptsize{x_{5}}$}
\put(39.7,12){$\scriptsize{x_{4}}$}
\put(45,11){$\scriptsize{x_{3}}$}
\put(52,11){$\scriptsize{x_{12}}$}
\put(58.5,13){$\scriptsize{x_{11}}$}
\put(65,10){$\scriptsize{x_{10}}$}

\put(25,-24){\scriptsize Figure 5 : $path(w_{i}\rightarrow w_{i+1}\rightarrow w_{i+2} \rightarrow w_{i+3})$}
\end{picture}

\vspace{6 cm}

\vspace{1 mm}
\begin{picture}(0,0)(-30,20)
\setlength{\unitlength}{1.5 mm}

\drawpolygon(0,0)(5,-5)(10,0)(10,5)(10,10)(5,15)(0,10)
\drawpolygon(10,0)(20,0)(20,10)(10,10)
\drawpolygon(20,0)(25,-5)(30,0)(30,5)(30,10)(25,15)(20,10)
\drawpolygon(30,0)(40,0)(40,10)(30,10)
\drawpolygon(40,0)(45,-5)(50,0)(50,5)(50,10)(45,15)(40,10)
\drawpolygon(50,0)(60,0)(60,10)(50,10)

\drawpolygon(0,0)(5,-5)(5,-15)(0,-10)
\drawpolygon(10,0)(20,0)(20,-10)(10,-10)
\drawpolygon(10,0)(5,-5)(5,-15)(10,-10)

\drawpolygon(20,0)(25,-5)(25,-15)(20,-10)
\drawpolygon(30,0)(40,0)(40,-10)(30,-10)
\drawpolygon(30,0)(25,-5)(25,-15)(30,-10)

\drawpolygon(40,0)(45,-5)(45,-15)(40,-10)
\drawpolygon(50,0)(60,0)(60,-10)(50,-10)
\drawpolygon(50,0)(45,-5)(45,-15)(50,-10)

\drawline[AHnb=0](5,-5)(0,-10)
\drawline[AHnb=0](5,-5)(10,-10)
\drawline[AHnb=0](10,0)(20,-10)
\drawline[AHnb=0](20,0)(10,10)

\drawline[AHnb=0](25,-5)(20,-10)
\drawline[AHnb=0](25,-5)(30,-10)
\drawline[AHnb=0](30,0)(40,-10)
\drawline[AHnb=0](40,0)(30,10)

\drawline[AHnb=0](45,-5)(40,-10)
\drawline[AHnb=0](45,-5)(50,-10)
\drawline[AHnb=0](50,0)(60,-10)
\drawline[AHnb=0](60,0)(50,10)

\put(-3,-10){$\scriptsize{v_{1}}$}
\put(5,-17){$\scriptsize{v_{2}}$}
\put(10,-12){$\scriptsize{v_{3}}$}
\put(18,-12){$\scriptsize{v_{4}}$}
\put(25,-17){$\scriptsize{v_{5}}$}
\put(30,-12){$\scriptsize{v_{6}}$}
\put(38,-12){$\scriptsize{v_{7}}$}
\put(45,-17){$\scriptsize{v_{8}}$}
\put(50,-12){$\scriptsize{v_{9}}$}
\put(60,-12){$\scriptsize{v_{10}}$}

\put(-4,0){$\scriptsize{u_{1}}$}
\put(4,-3){$\scriptsize{u_{2}}$}
\put(6,0){$\scriptsize{u_{3}}$}
\put(21,0){$\scriptsize{u_{4}}$}
\put(23,-3){$\scriptsize{u_{5}}$}
\put(26,0){$\scriptsize{u_{6}}$}
\put(41,0){$\scriptsize{u_{7}}$}
\put(44,-3){$\scriptsize{u_{8}}$}
\put(47,0){$\scriptsize{u_{9}}$}
\put(61,0){$\scriptsize{u_{10}}$}

\put(22,-24){\scriptsize Figure 6}

\end{picture}

\vspace{6 cm}

\begin{defn}\label{df4.2} Let $P(u_{1}\rightarrow u_{r})$ be a path in $M$. We denote three paths $Q_{1}(w_{i}\rightarrow w_{i+1}\rightarrow w_{i+2} \rightarrow w_{i+3})$, $Q_{2}(v_{i}\rightarrow v_{i+1}\rightarrow v_{i+2} \rightarrow v_{i+3})$ and $Q_{3}(z_{i}\rightarrow z_{i+1}\rightarrow z_{i+2} \rightarrow z_{i+3})$ where $link(w_{i})=C(w_{i-1},$ $x_{2},w_{i+1},$ $\textbf{w}_{i+2},$ $\textbf{x}_{3},$ $\textbf{x}_{4},$ $ x_{5},$ $x_{6})$, $link(w_{i+1})=C(w_{i},$ $x_{2},$ $x_{7},$ $x_{9},$ $w_{i+2},$ $\textbf{x}_{3},$ $ \textbf{x}_{4},$ $\textbf{x}_{5})$, $link(w_{i+2})=C(w_{i+1},$ $x_{8},$ $x_{9},$ $w_{i+3},$ $x_{3},$ $\textbf{x}_{4},$ $ \textbf{x}_{5},$ $\textbf{w}_{i},)$, $link(w_{i+3})=C(w_{i+2},$ $x_{9},$ $w_{i+4},$ $\textbf{w}_{i+5},$ $\textbf{x}_{10},$ $\textbf{x}_{11},$ $ x_{12},x_{3})$, $link(v_{i})=C(v_{i-1},$ $y_{1},y_{2},$ $y_{3},$ $v_{i+1},$ $\textbf{y}_{4},$ $ \textbf{y}_{5},$ $\textbf{y}_{6})$, $link(v_{i+1})=C(v_{i},$ $y_{3},$ $y_{7},$ $v_{i+2},$ $y_{4},$ $\textbf{y}_{5},$ $ \textbf{y}_{6},$ $\textbf{v}_{i-1})$, $link(v_{i+2})=C(v_{i+1},$ $v_{7},$ $v_{i+3},$ $\textbf{v}_{i+4},$ $\textbf{y}_{11},$ $\textbf{y}_{12},$ $ y_{8},$ $w_{4},)$, $link(v_{i+3})=C(v_{i+2},$ $y_{7},$ $y_{9},$ $y_{10},$ $v_{i+4},$ $\textbf{y}_{11},$ $ \textbf{y}_{12},\textbf{y}_{8})$, $link(z_{i})=C(z_{i-1},$ $r_{3},r_{4},$ $z_{i+1},$ $r_{12},$ $ \textbf{r}_{13},$ $\textbf{r}_{14}, \textbf{z}_{i-2})$, $link(z_{i+1})=C(z_{i},$ $r_{4},$ $z_{i+2},$ $\textbf{z}_{i+3},$ $\textbf{r}_{9},$ $\textbf{r}_{10},$ $ \textbf{r}_{11},$ $r_{12})$, $link(z_{i+2})=C(z_{i+1},$ $r_{4},$ $r_{5},$ $r_{6},$ $z_{i+3},$ $\textbf{r}_{9},$ $ \textbf{r}_{10},$ $\textbf{r}_{11},)$ and $link(z_{i+3})=C(z_{i+2},$ $r_{6},$ $r_{7},$ $z_{i+4},$ $r_{9},$ $\textbf{r}_{10},$ $ \textbf{r}_{11},\textbf{z}_{i+1})$. A path $P$ is said to be well defined if we consider a path of length three, say $L(u_{t}\rightarrow u_{t+1}\rightarrow u_{t+2}\rightarrow u_{t+3})$ in P or in the extended path of $P$ then either $L \mapsto Q_{1}$ by $u_{j}\mapsto w_{2t+3-j}$ or $L \mapsto Q_{2}$ by $u_{j}\mapsto v_{2t+3-j}$ or $L \mapsto Q_{3}$ by $u_{j}\mapsto z_{2t+3-j}$ for $j \in \{t,t+1,t+2,t+3\}$. In Figure 6, $P = path(v_{1}\rightarrow v_{2}\rightarrow \dots \rightarrow v_{10})$. We denote this path by $Y_{2}$.
\end{defn}

We consider maximal path of type $Y_{i}$ for $i \in  \{1, 2\}$. It is cycle by the similar argument of theorem \ref{thm1.1}. Let $C_{1}(u_{1}, u_{2}, \dots, u_{r})$ of type $Y_{1}$ and $C_{2}(v_{1}, v_{2}, \dots, v_{r})$ of type $Y_{2}$ be two cycles of same length $r$. Let $u_{i-1}\rightarrow u_{i} \rightarrow u_{i+1}$ be a sub path of $C_{1}$ where adjacent $6$-gon is on one side and all four $3$-gons are on the other side at $u_{i}$. Similarly, we consider a sub path $v_{j-1}\rightarrow v_{j} \rightarrow v_{j+1}$ of $C_{2}$ where adjacent $6$-gon is on one side and all four $3$-gons are on the other side at $v_{j}$. We define a map which maps $u_{i}$ to $v_{j}$. Then, by definition \ref{df4.1}, \ref{df4.2}, the vertex $u_{i+1}$ maps $v_{j-1}$ and $u_{i-1}$ maps $v_{j+1}$. We extend it to hole cycles. Hence we get $u_{i+t}\mapsto v_{j-t}$ for $0\leq t \leq r-1$. This map gives the cycles $C_{1}$ and $C_{2}$ are of type $D_{1} ( = D_{2})$. Therefore, the definition \ref{df4.1} and \ref{df4.2} defined same cycle.

\vspace{1 mm}
\begin{picture}(0,0)(-10,100)
\setlength{\unitlength}{1.5 mm}

\drawpolygon(0,0)(5,-5)(10,0)(10,5)(10,10)(5,15)(0,10)
\drawpolygon(10,0)(20,0)(20,10)(10,10)
\drawpolygon(20,0)(25,-5)(30,0)(30,5)(30,10)(25,15)(20,10)
\drawpolygon(30,0)(40,0)(40,10)(30,10)
\drawpolygon(40,0)(45,-5)(50,0)(50,5)(50,10)(45,15)(40,10)
\drawpolygon(50,0)(60,0)(60,10)(50,10)

\drawpolygon(0,0)(5,-5)(5,-15)(0,-10)
\drawpolygon(10,0)(20,0)(20,-10)(10,-10)
\drawpolygon(10,0)(5,-5)(5,-15)(10,-10)

\drawpolygon(20,0)(25,-5)(25,-15)(20,-10)
\drawpolygon(30,0)(40,0)(40,-10)(30,-10)
\drawpolygon(30,0)(25,-5)(25,-15)(30,-10)

\drawpolygon(40,0)(45,-5)(45,-15)(40,-10)
\drawpolygon(50,0)(60,0)(60,-10)(50,-10)
\drawpolygon(50,0)(45,-5)(45,-15)(50,-10)

\drawline[AHnb=0](5,-5)(0,-10)
\drawline[AHnb=0](5,-5)(10,-10)
\drawline[AHnb=0](10,0)(20,-10)
\drawline[AHnb=0](20,0)(10,10)

\drawline[AHnb=0](25,-5)(20,-10)
\drawline[AHnb=0](25,-5)(30,-10)
\drawline[AHnb=0](30,0)(40,-10)
\drawline[AHnb=0](40,0)(30,10)

\drawline[AHnb=0](45,-5)(40,-10)
\drawline[AHnb=0](45,-5)(50,-10)
\drawline[AHnb=0](50,0)(60,-10)
\drawline[AHnb=0](60,0)(50,10)

\put(-3,-10){$\scriptsize{v_{1}}$}
\put(5,-17){$\scriptsize{v_{2}}$}
\put(10,-12){$\scriptsize{v_{3}}$}
\put(18,-12){$\scriptsize{v_{4}}$}
\put(25,-17){$\scriptsize{v_{5}}$}
\put(30,-12){$\scriptsize{v_{6}}$}
\put(38,-12){$\scriptsize{v_{7}}$}
\put(45,-17){$\scriptsize{v_{8}}$}
\put(50,-12){$\scriptsize{v_{9}}$}
\put(60,-12){$\scriptsize{v_{1}}$}

\put(-4,0){$\scriptsize{v_{10}}$}
\put(4,-3){$\scriptsize{v_{11}}$}
\put(6,0){$\scriptsize{v_{12}}$}
\put(21,0){$\scriptsize{v_{13}}$}
\put(23,-3){$\scriptsize{v_{14}}$}
\put(26,0){$\scriptsize{v_{15}}$}
\put(41,0){$\scriptsize{v_{16}}$}
\put(44,-3){$\scriptsize{v_{17}}$}
\put(46,0){$\scriptsize{v_{18}}$}
\put(61,0){$\scriptsize{v_{10}}$}


\drawpolygon(5,25)(10,20)(15,25)(15,30)(15,35)(10,40)(5,35)
\drawpolygon(15,25)(25,25)(25,35)(15,35)
\drawpolygon(5,25)(10,20)(10,10)(5,15)
\drawpolygon(15,25)(25,25)(25,15)(20,10)
\drawpolygon(15,25)(10,20)(10,10)(20,10)
\drawline[AHnb=0](10,20)(5,15)
\drawline[AHnb=0](10,20)(20,10)
\drawline[AHnb=0](15,25)(25,15)
\drawline[AHnb=0](15,35)(25,25)

\drawpolygon(25,25)(30,20)(35,25)(35,30)(35,35)(30,40)(25,35)
\drawpolygon(35,25)(45,25)(45,35)(35,35)
\drawpolygon(25,25)(30,20)(30,10)(25,15)
\drawpolygon(35,25)(45,25)(45,15)(40,10)
\drawpolygon(35,25)(30,20)(30,10)(40,10)
\drawline[AHnb=0](30,20)(25,15)
\drawline[AHnb=0](30,20)(40,10)
\drawline[AHnb=0](35,25)(45,15)
\drawline[AHnb=0](35,35)(45,25)

\drawpolygon(45,25)(50,20)(55,25)(55,30)(55,35)(50,40)(45,35)
\drawpolygon(55,25)(65,25)(65,35)(55,35)
\drawpolygon(45,25)(50,20)(50,10)(45,15)
\drawpolygon(55,25)(65,25)(65,15)(60,10)
\drawpolygon(55,25)(50,20)(50,10)(60,10)
\drawline[AHnb=0](50,20)(45,15)
\drawline[AHnb=0](50,20)(60,10)
\drawline[AHnb=0](55,25)(65,15)
\drawline[AHnb=0](55,35)(65,25)

\put(1,15){$\scriptsize{v_{19}}$}
\put(7,8){$\scriptsize{v_{20}}$}
\put(16,8){$\scriptsize{v_{21}}$}
\put(23.5,12.5){$\scriptsize{v_{22}}$}
\put(26.5,8){$\scriptsize{v_{23}}$}
\put(40,9){$\scriptsize{v_{24}}$}
\put(44,12){$\scriptsize{v_{25}}$}
\put(47,8){$\scriptsize{v_{26}}$}
\put(61,9){$\scriptsize{v_{27}}$}
\put(64.8,13){$\scriptsize{v_{19}}$}

\put(4,5){$\scriptsize F_{1}$}
\put(24,5){$\scriptsize F_{2}$}
\put(44,5){$\scriptsize F_{3}$}

\put(9,30){$\scriptsize F_{4}$}
\put(29,30){$\scriptsize F_{5}$}
\put(49,30){$\scriptsize F_{6}$}

\put(14,55){$\scriptsize F_{7}$}
\put(34,55){$\scriptsize F_{8}$}
\put(54,55){$\scriptsize F_{9}$}

\put(1,25){$\scriptsize{v_{28}}$}
\put(9,22){$\scriptsize{v_{29}}$}
\put(11,25){$\scriptsize{v_{30}}$}
\put(26,25){$\scriptsize{v_{31}}$}
\put(28,22){$\scriptsize{v_{32}}$}
\put(31,25){$\scriptsize{v_{33}}$}
\put(46,25){$\scriptsize{v_{34}}$}
\put(49,22){$\scriptsize{v_{35}}$}
\put(51,25){$\scriptsize{v_{36}}$}
\put(66,25){$\scriptsize{v_{28}}$}

\put(-3,10){$\scriptsize{v_{27}}$}


\drawpolygon(10,50)(15,45)(20,50)(20,55)(20,60)(15,65)(10,60)
\drawpolygon(20,50)(30,50)(30,60)(20,60)
\drawpolygon(10,50)(15,45)(15,35)(10,40)
\drawpolygon(20,50)(30,50)(30,40)(25,35)
\drawpolygon(20,50)(15,45)(15,35)(25,35)
\drawline[AHnb=0](15,45)(10,40)
\drawline[AHnb=0](15,45)(25,35)
\drawline[AHnb=0](20,50)(30,40)
\drawline[AHnb=0](20,60)(30,50)

\drawpolygon(30,50)(35,45)(40,50)(40,55)(40,60)(35,65)(30,60)
\drawpolygon(40,50)(50,50)(50,60)(40,60)
\drawpolygon(30,50)(35,45)(35,35)(30,40)
\drawpolygon(40,50)(50,50)(50,40)(45,35)
\drawpolygon(40,50)(35,45)(35,35)(45,35)
\drawline[AHnb=0](35,45)(30,40)
\drawline[AHnb=0](35,45)(45,35)
\drawline[AHnb=0](40,50)(50,40)
\drawline[AHnb=0](40,60)(50,50)

\drawpolygon(50,50)(55,45)(60,50)(60,55)(60,60)(55,65)(50,60)
\drawpolygon(60,50)(70,50)(70,60)(60,60)
\drawpolygon(50,50)(55,45)(55,35)(50,40)
\drawpolygon(60,50)(70,50)(70,40)(65,35)
\drawpolygon(60,50)(55,45)(55,35)(65,35)
\drawline[AHnb=0](55,45)(50,40)
\drawline[AHnb=0](55,45)(65,35)
\drawline[AHnb=0](60,50)(70,40)
\drawline[AHnb=0](60,60)(70,50)

\put(6,40){$\scriptsize{v_{37}}$}
\put(11,34){$\scriptsize{v_{38}}$}
\put(25,33){$\scriptsize{v_{39}}$}
\put(29,37){$\scriptsize{v_{40}}$}
\put(31.5,33){$\scriptsize{v_{41}}$}
\put(45,33){$\scriptsize{v_{42}}$}
\put(49,37){$\scriptsize{v_{43}}$}
\put(51.5,33){$\scriptsize{v_{44}}$}
\put(65,33){$\scriptsize{v_{45}}$}
\put(70,38){$\scriptsize{v_{37}}$}

\put(6,50){$\scriptsize{v_{46}}$}
\put(14,47){$\scriptsize{v_{47}}$}
\put(16,50){$\scriptsize{v_{48}}$}
\put(31,50){$\scriptsize{v_{49}}$}
\put(33,47){$\scriptsize{v_{50}}$}
\put(36,50){$\scriptsize{v_{51}}$}
\put(51,50){$\scriptsize{v_{52}}$}
\put(54,47){$\scriptsize{v_{53}}$}
\put(56,50){$\scriptsize{v_{54}}$}
\put(71,50){$\scriptsize{v_{46}}$}

\put(2,35){$\scriptsize{v_{45}}$}

\put(6,60){$\scriptsize{v_{3}}$}
\put(14,66){$\scriptsize{v_{4}}$}
\put(16,60){$\scriptsize{v_{5}}$}
\put(31,60){$\scriptsize{v_{6}}$}
\put(33,66){$\scriptsize{v_{7}}$}
\put(36,60){$\scriptsize{v_{8}}$}
\put(51,60){$\scriptsize{v_{9}}$}
\put(54,66){$\scriptsize{v_{1}}$}
\put(57,60){$\scriptsize{v_{2}}$}
\put(71,60){$\scriptsize{v_{3}}$}

\put(20,-25){\scriptsize Figure 7 : $T(9, 6, 2)$}

\end{picture}

\vspace{15 cm}

Let $M$ be a map. We construct a cylinder using a cycle of type $Y_{1}$ like done in Section 2. We get a cylinder $S_{C_{i}}$ by a cycle $C_{i}$. It's boundary cycles are, say $C_{l}, C_{k}$. By the similar argument of lemma \ref{lem1.2}, $length(C_{i}) = length(C_{l}) = length(C_{k})$. Also, by the similar argument of lemma \ref{lem1.3}, the cycles of type $Y_{1}$ which are homologous to $C_{1}$ have same length. Its clear from the definition and representation that there are three cycles at each vertex, say $v$, of type $Y_{1}$. Let $L_{1}(v), L_{2}(v), L_{3}(v)$ be three cycles at $v$. We define $T(r, s, k)$ of $M$ by the similar argument of Section 2. In this process, we take second cut along $L_{2}$ in place of $L_{3}$ where the starting adjacent face to horizontal base cycle $L_{1}$ is $3$-gon. Hence we get a $(r, s, k)$-representation of $M$. The Figure 7 is an example of $T(9, 6, 3)$ of a map of type $\{ 3^{4}, 6\}$ on 54 vertices. In Figure 7, $L_{1} = C(v_{1},v_{2},\dots,v_{9})$ and the path $v_{1}\rightarrow v_{10}\rightarrow v_{27}\rightarrow v_{19}\rightarrow v_{28}\rightarrow v_{45}\rightarrow v_{37}\rightarrow v_{46}\rightarrow v_{3}$ is part of the cycle $L_{2}$.  This $T(r, s, k)$ exists for every $M$ of type $\{ 3^{4}, 6\}$. Now, we show the map $M$ contains contractible Hamiltonian by theorem \ref{thm4.1}. Therefore the map of type $\{ 3^{4}, 6\}$ on the torus is Hamiltonian.


\begin{theo} \label{thm4.1} Let $M$ be a map of type $\{ 3^{4}, 6\}$ on the torus. Then it contains a contractible Hamiltonian cycle.
\end{theo}

\noindent{\sc Proof of Theorem}\ref{thm4.1} Let $M$ be a map of type $\{ 3^{4}, 6\}$ on the torus. We consider a vertex, say $v$. Its clear from type of $M$ that the vertex $v$ belongs to exactly one $6$-gon, say $F_{v,6}$. In Figure 7, if we take any $F_{i}$ and $F_{j}$ for $i\not=j$ then $F_{i}$, $F_{i}$ are disjoint. We choose any two $6$-gons, say $F_{u,6}$, $F_{v,6}$ in $M$. The faces $F_{u,6}$ and $F_{v,6}$ are disjoint. We consider all $6$-gons of $M$, say $F_{v_{1},6}, F_{v_{2},6},\dots, F_{v_{t},6}$ where $V(F_{v_{i},6})\cap V(F_{v_{j},6})=\emptyset$ for $i \not= j$. That is, $\cup_{i=1}^{t}V(F_{v_{i},6})=V(M)$. Let $C_{1}, C_{2},\dots,C_{t}$ be a list of homologous horizontal cycles of type $Y_{1}$ in $T(r, s, k)$. In Figure 7, $s = 6$ and $C_{1} = C(v_{1},$ $v_{2},$ $\dots,$ $v_{9})$, $C_{2} = C(v_{10},$ $v_{11},$ $\dots,$ $v_{18})$, $C_{3} = C(v_{19},$ $v_{20},$ $\dots,$ $v_{27})$, $C_{4} = C(v_{28},$ $v_{29},$ $\dots,$ $v_{36})$, $C_{5} = C(v_{37},$ $v_{38},$ $\dots,$ $v_{45})$, $C_{6} = C(v_{46},$ $v_{47},$ $\dots,$ $v_{54})$. For $i, j \in \{1,2,\dots,t\}$ and $|i - j| \not= 1$, we consider two cycles $C_{i}, C_{j}$. We denote $W_{C,6}$ be a set of adjacent $6$-gons to $C$. In Figure 7, $W_{C_{2},6} = \{F_{1}, F_{2}, F_{3}\}$. Then its clear from the representation, $W_{C_{i},6}\cap W_{C_{j},6} = \emptyset$. We define a set $S_{C_{i}} \colon=\{ \phi \}$. We choose $C_{i}$ and consider all adjacent $6$-gons in $S_{C_{i}}$. Let $F_{u,6}$, $F_{v,6}$ be two consecutive $6$-gons adjacent to $C_{i}$. Its clear from the representation $T(r, s, k)$, there are two $3$-gons, say $\triangle_{1}, \triangle_{2}$ such that geometric carrier of $F_{u, 6}, F_{v, 6}, \triangle_{1}, \triangle_{2}$ is $2$-disk, $|\{\triangle_{1}, \triangle_{2}\}|\cap F_{u,6}$ is an edge, $|\{\triangle_{1}, \triangle_{2}\}|\cap F_{v,6}$ is  an edge and $|\{\triangle_{1}, \triangle_{2}\}|\cap F_{w,6} = \emptyset$ for all $w \in V(M)\setminus \{u,v\}$. In Figure 7, if we consider $F_{1}$ and $F_{2}$ then triangles are $\triangle_{1}(v_{12},v_{13},v_{20})$, $\triangle_{2}(v_{13},v_{21},v_{20})$. We put these two triangles $\triangle_{1}, \triangle_{2}$ in $S_{c_{i}}$. We do it for any two consecutive $6$-gons which are adjacent to $C_{i}$ except first and last one. In Figure 7, if we consider $C_{2}$ then first and last faces are $F_{1}$ and $F_{3}$ respectively. Hence we get a set $S_{C_{i}}$ of $6$-gons and $3$-gons where the geometric carrier is a $2$-disk and bounded by a cycle. Therefore we get a set $S_{C_{i}}$ for cycle $C_{i}$ where $i$ is even and $i \in \{2, 4, 6, \dots , t\}$. We consider the set $S := \cup_{i=2}^{t}S_{C_{i}}$. The set $S$ contains all $6$-gons of $M$. Now, we choose a vertex of $C_{i_{0}}$ for a fixed $i_{0}$ where $1 \leq i_{0} \leq t$. We consider a path, say $Q$ of type $Y_{1}$ vertically which is not homologous to $C_{i_{0}}$. Similarly, we choose two adjacent $3$-gons on one side of $Q$, say $\triangle_{i}, \triangle_{i+2}$, between any two consecutive $6$-gons adjacent to $C_{i}$ and $C_{i+2}$ along $Q$. In Figure 7, $Q = path(v_{1}\rightarrow v_{10} \rightarrow v_{27} \rightarrow v_{19} \rightarrow v_{28} \rightarrow v_{45} \rightarrow v_{37} \rightarrow v_{46} \rightarrow v_{3})$. So, we get $|S_{C_{i}} \cup S_{C_{i+2}}\cup \{\triangle_{i}, \triangle_{i+2}\}|$ which is $2$-disk. We collect all triangle between any two consecutive $6$-gons along $Q$ in $S$. Hence the geometric carrier $|S|$ which is a $2$-disk and bounded by a cycle. The set $S$ contains $F_{v_{1},6}, F_{v_{2},6},\dots, F_{v_{t},6}$ and $|S|$ is $2$-disk. Therefore $\partial |S|$ is a cycle and contains all the vertices of $M$. Hence the map $M$ contains contractible Hamiltonian cycle.
\hfill$\Box$

\section{Map of type $\{ 4, 8, 8\}$}

 Let $M$ be a map of type $\{ 4, 8, 8\}$. We denote $C_{F_{l}}$ to be a cycle of length $l$ which is boundary of a face $F_{l}$. We use this notation in the definition \ref{df5.1}, \ref{df6.1}, \ref{df7.1} and \ref{df8.1}. We define a path in $M$ as follows $\colon$

\begin{defn}\label{df5.1} Let $P(u_{1}\rightarrow u_{r})$ be a path in $M$. The path $P$ is said to well defined if it satisfy the following properties $\colon$
\begin{enumerate}

\item if $P \cap C_{F_{4}} \neq \emptyset$ then $P\cap C_{F_{4}}$ is a path of length one,

\item if $P\cap C_{F_{8}} \neq \emptyset$ then $P\cap C_{F_{8}}$ is a path of length three.

\item For the boundary vertices, we consider the extended path of $P$.

\end{enumerate}
We denote this path by $Z_{1}$.
\end{defn}

We consider a maximal path of type $Z_{1}$ which is cycle by the similar argument of theorem \ref{thm1.1}. We define a cylinder with two boundaries where boundaries are cycles of type $Z_{1}$ by the construction done in Section 5. Let $C_{i}$ be a cycle. Let $S_{C_{i}}$ be a cylinder by $C_{i}$. We consider the boundaries of $S_{C_{i}}$, say $C_{l}, C_{k}$. By the similar argument of lemma \ref{lem1.2}, $length(C_{i}) = length(C_{l}) = length(C_{k})$. Let $C$ be a cycle of type $Z_{1}$ which is homologous to $C_{1}$. Then, by the similar argument of lemma \ref{lem1.3}, $length(C) = length(C_{1})$. Its clear from the definition \ref{df5.1} that there are two cycles of type $Z_{1}$ at a vertex. Therefore, we consider two cycles of type $Z_{1}$ at a vertex $v$, say $L_{1}(v), L_{2}(v)$. We define a $(r, s, k)$-representation of $M$ by similar idea of Section 3. Here we take second cut along cycle $L_{2}$ where the starting adjacent face to the base horizontal cycle $L_{1}$ is $4$-gon. Hence we get a $(r, s, k)$-representation of $M$. Therefore, every map of type $\{ 4, 8, 8\}$ on the torus has $T(r, s, k)$. Now the map $M$ has contractible Hamiltonian cycle by theorem \ref{thm5.1}. Therefore any map of type $\{ 4, 8, 8\}$ on the torus is Hamiltonian.

\begin{theo} \label{thm5.1} Let $M$ be a map of type $\{ 4, 8, 8\}$ on the torus. Then it contains a contractible Hamiltonian cycle.
\end{theo}

\noindent{\sc Proof of Theorem}\ref{thm5.1} Let $M$ be a map. Let $C_{1}, C_{2}, \dots, C_{s}$ be the list of all homologous horizontal cycles of type $Z_{1}$ in $T(r,s,k)$ of length, say $r$. Let $C_{1}(u_{1_{1}}, u_{1_{2}}, \dots, u_{1_{r}})$, $C_{2}(u_{2_{1}},$ $u_{2_{2}},$ $\dots, u_{2_{r}})$, $C_{3}(u_{3_{1}},$ $ u_{3_{2}},$ $\dots, u_{3_{r}})$, \dots, $C_{s}(u_{s_{1}}, u_{s_{2}},\dots u_{s_{r}})$ be the cycles in order, that is, if we consider any consecutive three cycles $C_{k}, C_{k+1}$ and $C_{k+2}$ then $S_{C_{k+1}}$ is a cylinder with $\partial S_{C_{k+1}} = \{ C_{k}, C_{k+2}\}$. We define a cycle $C_{i-1,i} \colon = C(u_{i-1_{2}},$ $\dots,$ $u_{i-1_{r}},$ $u_{i-1_{1}},$ $u_{i_{1}},$ $u_{i_{r}},$ $u_{i_{r-1}},$ $\dots,$ $u_{i_{2}})$ by $C_{i-1}$ and $C_{i}$. The cycle $C_{i-1,i}$ is contractible as it bounds a $2$-disk which consists of $4$-gons and $8$-gons. We follow one of the following cases $\colon$

\begin{enumerate}
\item Let $s$ be an even integer and the cycle $C_{i-1, i}$ for $i = 2, 4,\dots, s$. We consider two cycles $C_{i-1, i}$, $C_{i+1, i+2}$ where $i$ is even and adjacent $4$-gon $C_{F_{4}}(u_{i_{3}},u_{i_{4}}, u_{i+1_{4}},u_{i+1_{3}})$. Here $C_{F_{4}}\cap C_{i-1, i}$ and $C_{F_{4}}\cap C_{i+1, i+2}$ are two edges. Therefore we define a cycle $C_{i-1,i}\cup C_{i+1,i+2} := C(u_{i-1_{2}},$ $\dots,$ $u_{i-1_{r}},$ $u_{i-1_{1}},$ $u_{i_{1}},$ $u_{i_{r}},$ $u_{i_{r-1}},$ $\dots, u_{i_{4}}$ $u_{i+1_{4}},$ $u_{i+1_{5}},$ $\dots,$ $u_{i+1_{r}},$ $u_{i+1_{1}},$ $u_{i+2_{1}},$ $u_{i+2_{r}},$ $u_{i+2_{r-1}},$ $\dots$, $u_{i+2_{2}},$ $u_{i+1_{2}},$ $u_{i+1_{3}},$ $u_{i_{3}},$ $u_{i_{2}})$. By this definition we define cycle $C := C_{1,2}\cup C_{3,4} \cup C_{5,6}\cup \dots \cup C_{s-2,s}$ which is contractible. It contains all the vertices of $C_{1}, C_{2}, \dots, C_{s}$. Therefore $M$ contains a contractible Hamiltonian cycle.

\item Let $s$ be an odd integer. When $s=1$ we define $C := C_{1}$. When $s > 1$, we have a contractible cycle $C_{i-1,i}$ for $i = 2, 4, \dots, s-1$. Next we use the idea $(1)$. We define a cycle $C := C_{1,2}\cup C_{3,4} \cup C_{5,6}\cup \dots \cup C_{s-2,s-1}$. Let $F_{4}$ be a $4$-gon adjacent along an edge with $C_{s}$ and $C_{s-1}$. We concatenate the cycle $C$ with $C_{F_{4}}$. We do it for all adjacent $4$-gons those are adjacent with $C_{s}$ and $C_{s-1}$. This gives a cycle which is contractible and contains all the vertices of $M$. Therefore the map $M$ contains a contractible Hamiltonian cycle.
\hfill$\Box$
\end{enumerate}

\section{Map of type $\{ 3, 12, 12\}$}

 Let $M$ be a map of type $\{ 3, 12, 12\}$. We define a path as follows $\colon$

\begin{defn}\label{df6.1} Let $P(u_{1}\rightarrow u_{r})$ be a path. The path $P$ is said to well defined if it follow the following properties $\colon$
\begin{enumerate}

\item if $P \cap C_{F_{3}} \neq \emptyset$ then $P \cap C_{F_{3}}$ is a path of length one,

\item if $P \cap C_{F_{12}} \neq \emptyset$ then $P \cap C_{F_{12}}$ is a path of length three.

\item  For the boundary vertices, we consider extended path of $P$.

\end{enumerate}
We denote this path by $G_{1}$.
\end{defn}

We consider a maximal path of type $G_{1}$ which is cycle by similar argument of theorem \ref{thm1.1}. Therefore any maximal path of type $G_{1}$ is cycle. Its clear from the definition \ref{df6.1} that there are two cycles of type $G_{1}$ at a vertex. Let $L_{1}(v), L_{2}(w), L_{3}(v)$ be three cycles at $v,$ $w$ where $vw$ is an edge of $L_{1}$ and belongs to a $3$-gon. We define a $(r, s, k)$-representation of $M$ by similar idea of Section 4. Here we take second cut along a cycle of type $G_{1}$ where the starting adjacent face to horizontal base cycle $L_{1}$ is $12$-gon. Therefore we get $(r, s, k)$-representation of $M$. We give an example \ref{example2} which is not Hamiltonian by showing that it does not contain any maximal contractible or non-contractible cycle which contains all the vertices of the map $M$. Therefore the map of type $\{ 3, 12, 12\}$ on the torus is not Hamiltonian.

\begin{example}\label{example2}
Example of a map of type $\{3, 12, 12\}$ on the torus which is non-Hamiltonian.
\end{example}

Let $M$ be a map which has a representation $T(24, 2, 9)$. We consider a maximal non contractible cycle of type $G_{1}$ in $T(24, 2, 9)$. It contains maximal cycle of length $24$ among all possible cycles of type $G_{1}$. Let $C(1, 2,\dots, 24)$ be a cycle of length $24$. Now we extend $C$ by concatenate a cycle of length $l$ which bounds a $l$-gon. The representation $T(24, 2, 9)$ contains only $12$-gons and $3$-gons. Therefore we can extend $C$ to a cycle of length $24 + 6 = 30$ for a $12$-gon. Similarly, we can extend it to a cycle of length $24 + 1 = 25$ for a $3$-gon. Therefore, the length of the extended cycle of $C$ is, say $ 24 + 6 m_{1} + m_{2}$ where $m_{1}$ denotes the number of $12$-gons and $m_{2}$ denotes the number of $3$-gons. In Figure 8, we use the notation $F_{i, 12}$ to denote $i^{th}$ $12$-gon. Its clear from Figure 8, we can concatenate the cycles $\partial F_{1,12}$, $\partial F_{3,12}$ and $\partial F_{5,12}$ with $C$. If we consider any other $F_{i,12}$ for $i \not= 1, 3, 5$ then there will be a vertex which is repetition in the extended cycle. Therefore $m_{1} = 3$. Similarly, we have $m_{2} = 15 $ for $3$-gons. Hence the value of $24+6m_{1}+m_{2}$ is at most $57$ as $(24+6.3+15) = (24 + 18 + 15)=57 < 72$. Similarly, if we consider only $3$-gons or $12$-gons then also we get a cycle of length less than $72$.

Next, we consider a non contractible cycle which does not follow the definition {\ref{df6.1}. We extend this cycle to a cycle of maximal length by concatenating its boundary cycles of adjacent faces. We consider longest cycle among these maximal cycles. We denote it by $C$. In Figure 8, $C :=C(1,$ $24,$ $25,$ $31,$ $37,$ $38,$ $39,$ $40,$ $62,$ $68,$ $15,$ $14,$ $13,$ $12,$ $28,$ $34,$ $50,$ $49,$ $48,$ $47,$ $64,$ $70,$ $23,$ $22,$ $21,$ $20,$ $30,$ $ 36,$ $57,$ $58,$ $59,$ $60,$ $61,$ $67,$ $11,$ $10,$ $9,$ $27,$ $8,$ $7,$ $72,$ $6,$ $5,$ $4,$ $26,$ $32,$ $41,$ $42,$ $43,$ $44,$ $63,$ $69,$ $19,$ $18,$ $17,$ $416,$ $29,$ $35,$ $54,$ $53,$ $52,$ $51,$ $65,$ $ 71,$ $3,$ $2)$ of length $66 ( < 72)$. Therefore $T(24, 2, 9)$ does not have non contractible Hamiltonian cycle.

We choose a largest contractible cycle in $T(24, 2, 9)$. We denote $D_{C}$ to be a $2$-disk which is bounded by $C$. Let $2$-disk $D_{C_{1}}$ contains $m_{1}$ number of $12$-gons and $m_{2}$ number of $3$-gons. Then the length of the cycle $C_{1}$ is $12 + 10(m_{1}-1) + m_{2} = 10m_{1} + m_{2} + 2$. Now the number of $12$-gons can be at most $5$, that is, $m_{1} \leq 5$ otherwise we get a repetition of vertex in $C_{1}$. Similarly, the number of $3$-gons can be at most $14$, that is, $m_{2} \leq 14$ as we have at most $14$ adjacent $3$-gons to any $2$-disk which contains five $12$-gons. If we consider $m_{1} > 5$ or $m_{2} > 14$ then we get some vertex which are either inside the $2$-disk or repetition in $C_{1}$. Therefore, the value of $10m_{1}+m_{2}+2$ is at most $66$ as $(10m_{1}+m_{2}+2) \leq (50 + 14 + 2) = 66 < 72$. Hence $T(24, 2, 9)$ does not contain any contractible cycle which is Hamiltonian. Therefore the map $T(24, 2, 9)$ is not Hamiltonian.

\vspace{2 cm}
\smallskip
\begin{picture}(50,-40)(-2,0)
\setlength{\unitlength}{4.5mm}
\drawpolygon(-0.5,0)(1,-1)(2,-1)(3.5,0)(4,1)(4,2)(3.3,3)(2,4)(1,4)(-0.5,3)(-1,2)(-1,1)
\drawpolygon(4.5,0)(6,-1)(7,-1)(8.5,0)(9,1)(9,2)(8.3,3)(7,4)(6,4)(4.5,3)(4,2)(4,1)
\drawpolygon(9.5,0)(11,-1)(12,-1)(13.5,0)(14,1)(14,2)(13.3,3)(12,4)(11,4)(9.5,3)(9,2)(9,1)
\drawpolygon(14.5,0)(16,-1)(17,-1)(18.5,0)(19,1)(19,2)(18.3,3)(17,4)(16,4)(14.5,3)(14,2)(14,1)
\drawpolygon(19.5,0)(21,-1)(22,-1)(23.5,0)(24,1)(24,2)(23.3,3)(22,4)(21,4)(19.5,3)(19,2)(19,1)
\drawpolygon(24.5,0)(26,-1)(27,-1)(28.5,0)(29,1)(29,2)(28.3,3)(27,4)(26,4)(24.5,3)(24,2)(24,1)

\drawpolygon (-1,2)(-1.7,3)(-0.5,3)
\drawpolygon (4,2)(3.3,3)(4.5,3)
\drawpolygon (9,2)(8.3,3)(9.5,3)
\drawpolygon (14,2)(13.3,3)(14.5,3)
\drawpolygon (19,2)(18.3,3)(19.5,3)
\drawpolygon (24,2)(23.3,3)(24.5,3)

\drawpolygon (1,-1)(2,-1)(1.5,-2)
\drawpolygon (1.5,-2)(2,-1)(3.5,0)(4.5,0)(6,-1)(6.5,-2)(6.5,-3)(5.9,-4)(4.7,-5)(3,-5)(1.8,-4)(1.5,-3)
\drawpolygon (6.5,-2)(7,-1)(8.5,0)(9.5,0)(11,-1)(11.5,-2)(11.5,-3)(10.9,-4)(9.7,-5)(8,-5)(6.8,-4)(6.5,-3)
\drawpolygon (11.5,-2)(12,-1)(13.5,0)(14.5,0)(16,-1)(16.5,-2)(16.5,-3)(15.9,-4)(14.7,-5)(13,-5)(11.8,-4)(11.5,-3)
\drawpolygon (16.5,-2)(17,-1)(18.5,0)(19.5,0)(21,-1)(21.5,-2)(21.5,-3)(20.9,-4)(19.7,-5)(18,-5)(16.8,-4)(16.5,-3)
\drawpolygon (21.5,-2)(22,-1)(23.5,0)(24.5,0)(26,-1)(26.5,-2)(26.5,-3)(25.9,-4)(24.7,-5)(23,-5)(21.8,-4)(21.5,-3)
\drawpolygon (26.5,-2)(27,-1)(28.5,0)(29.5,0)(31,-1)(31.5,-2)(31.5,-3)(30.9,-4)(29.7,-5)(28,-5)(26.8,-4)(26.5,-3)
\drawpolygon (28.5,0)(29.5,0)(29,1)

\drawpolygon (31.5,-3)(30.9,-4)(32,-4)
\drawpolygon (26.5,-3)(25.9,-4)(26.8,-4)
\drawpolygon (21.5,-3)(20.9,-4)(21.8,-4)
\drawpolygon (16.5,-3)(15.9,-4)(16.8,-4)
\drawpolygon (11.5,-3)(10.9,-4)(11.8,-4)
\drawpolygon (6.5,-3)(5.9,-4)(6.8,-4)

\put(1.7,-2){\scriptsize 31}
\put(2.8,0){\scriptsize 39}
\put(4.5,0){\scriptsize 40}
\put(5.9,-.8){\scriptsize 41}
\put(6.5,-3){\scriptsize 26}
\put(5.6,-4){\scriptsize 4}
\put(4.7,-4.7){\scriptsize 3}
\put(3,-4.9){\scriptsize 2}
\put(1.8,-3.8){\scriptsize 1}
\put(1.5,-3){\scriptsize 25}

\put(6.7,-2){\scriptsize 32}
\put(6.9,-.7){\scriptsize 42}
\put(7.9,0){\scriptsize 43}
\put(9.5,0){\scriptsize 44}
\put(10.8,-.9){\scriptsize 45}
\put(11.5,-3){\scriptsize 27}
\put(10.5,-4){\scriptsize 8}
\put(9.7,-4.6){\scriptsize 7}
\put(8,-4.9){\scriptsize 6}
\put(6.8,-4){\scriptsize 5}

\put(11.7,-2){\scriptsize 33}
\put(11.6,-.8){\scriptsize 46}
\put(12.7,0){\scriptsize 47}
\put(14.5,0){\scriptsize 48}
\put(16,-.9){\scriptsize 49}
\put(16.5,-3){\scriptsize 28}
\put(15.2,-4){\scriptsize 12}
\put(14.7,-4.6){\scriptsize 11}
\put(13,-5){\scriptsize 10}
\put(11.8,-4){\scriptsize 9}

\put(16.7,-2){\scriptsize 34}
\put(16.8,-.8){\scriptsize 50}
\put(17.8,0){\scriptsize 51}
\put(19.5,0){\scriptsize 52}
\put(20.5,-.7){\scriptsize 53}
\put(21.5,-3){\scriptsize 29}
\put(20.2,-4){\scriptsize 16}
\put(19.6,-4.6){\scriptsize 15}
\put(18,-4.9){\scriptsize 14}
\put(16.8,-3.9){\scriptsize 13}

\put(21.7,-2){\scriptsize 35}
\put(21.8,-.6){\scriptsize 54}
\put(22.8,0){\scriptsize 55}
\put(24.5,0){\scriptsize 56}
\put(26,-.8){\scriptsize 57}
\put(26.5,-3){\scriptsize 30}
\put(25.3,-4){\scriptsize 20}
\put(24.4,-4.7){\scriptsize 19}
\put(23,-4.9){\scriptsize 18}
\put(21.8,-4){\scriptsize 17}

\put(26.7,-2){\scriptsize 36}
\put(26.8,-.7){\scriptsize 58}
\put(27.8,0){\scriptsize 59}
\put(29.5,0){\scriptsize 60}
\put(31,-1){\scriptsize 37}
\put(31.5,-2){\scriptsize 31}
\put(31.5,-3){\scriptsize 25}
\put(30.2,-4){\scriptsize 24}
\put(29.5,-4.7){\scriptsize 23}
\put(28,-4.7){\scriptsize 22}
\put(26.8,-4){\scriptsize 21}

\put(32,-4){\scriptsize 1}

\put(-0.5,0){\scriptsize 60}
\put(.7,-.8){\scriptsize 37}
\put(1.8,-.7){\scriptsize 38}
\put(4,1){\scriptsize 62}
\put(4.2,2){\scriptsize 68}
\put(3,3.1){\scriptsize 14}
\put(2,4){\scriptsize 13}
\put(1,4.1){\scriptsize 12}
\put(-1,3.1){\scriptsize 11}
\put(-.8,2){\scriptsize 67}
\put(-1,1){\scriptsize 61}

\put(-2,3.1){\scriptsize 10}

\put(9.2,2){\scriptsize 69}
\put(8,3.1){\scriptsize 18}
\put(7,4){\scriptsize 17}
\put(6,4.1){\scriptsize 16}
\put(4,3.1){\scriptsize 15}

\put(13,3.1){\scriptsize 22}
\put(12,4){\scriptsize 21}
\put(11,4.1){\scriptsize 20}
\put(9,3.1){\scriptsize 19}
\put(9,1){\scriptsize 63}

\put(18.3,3.1){\scriptsize 2}
\put(17,4){\scriptsize 1}
\put(16,4.1){\scriptsize 24}
\put(14,3.1){\scriptsize 23}
\put(14.2,2){\scriptsize 70}
\put(14,1){\scriptsize 64}

\put(23.3,3.1){\scriptsize 6}
\put(22,4){\scriptsize 5}
\put(21,4.1){\scriptsize 4}
\put(19,3.1){\scriptsize 3}
\put(19.2,2){\scriptsize 71}
\put(19,1){\scriptsize 65}

\put(29,1){\scriptsize 61}
\put(29,2){\scriptsize 67}
\put(28,3.2){\scriptsize 10}
\put(27,4){\scriptsize 9}
\put(26,4.1){\scriptsize 8}
\put(24.5,3.1){\scriptsize 7}
\put(24.2,2){\scriptsize 72}
\put(24,1){\scriptsize 66}

\put(3.1,-2.2){\scriptsize $F_{1,12}$}
\put(8.8,-2.2){\scriptsize $F_{2,12}$}
\put(13.8,-2.2){\scriptsize $F_{3,12}$}
\put(18.9,-2.2){\scriptsize $F_{4,12}$}
\put(23.9,-2.2){\scriptsize $F_{5,12}$}
\put(28.1,-2.2){\scriptsize $F_{6,12}$}

\put(1.1,1){\scriptsize $F_{7,12}$}
\put(6.1,1){\scriptsize $F_{8,12}$}
\put(11.1,1){\scriptsize $F_{9,12}$}
\put(16.1,1){\scriptsize $F_{10,12}$}
\put(21.1,1){\scriptsize $F_{11,12}$}
\put(26.1,1){\scriptsize $F_{12,12}$}

\put(10,-7){\scriptsize  Figure 8 : $T(24, 2, 9)$}

\end{picture}
\vspace{3 cm}

\begin{theo} \label{thm6.1} Let $M$ be a map of type $\{3, 12, 12\}$ on the torus. Then it is not Hamiltonian.
\end{theo}

\noindent{\sc Proof of Theorem}\ref{thm6.1} The proof of the theorem \ref{thm6.1} follows from the example \ref{example2}.
\hfill$\Box$

\section{Map of type $\{ 4, 6, 12\}$}

Let $M$ be a map of type $\{ 4, 6, 12\}$. We define a path as follows $\colon$

\begin{defn}\label{df7.1} Let $P(u_{1}\rightarrow u_{r})$ be a path in $M$. The path $P$ is said to be well defined if it follow the following properties $\colon$
\begin{enumerate}

\item if $P \cap C_{F_{4}} \neq \emptyset$ then $P\cap C_{F_{4}}$ is a path of length one,

\item if $P \cap C_{F_{6}} \neq \emptyset$ then $P\cap C_{F_{6}}$ is a path of length three,

\item if $P\cap C_{F_{12}} \neq \emptyset$ then $P\cap C_{F_{12}}$ is a path of length five.

\item the boundary vertex of $P$ also well defined in extended path of $P$.

\end{enumerate}
We denote this path by $H_{1}$.
\end{defn}

We consider a maximal path of type $H_{1}$. The path is cycle by similar argument of theorem \ref{thm1.1}.  We define a cylinder with two boundaries where boundaries are cycles of type $H_{1}$ by the construction done in Section 7. We denote it by $S_{C_{i}}$ for a cycle $C_{i}$. Let $\partial S_{C_{i}} = \{C_{l}, C_{k}\}$ where $C_{l}, C_{k}$ are of type $H_{1}$. By the lemma \ref{lem1.2}, $length(C_{i}) = length(C_{l}) = length(C_{k})$. Now we consider all possible horizontal homologous cycles $C_{1}, C_{2},\dots,C_{m}$ of type $H_{1}$ in $M$. By the similar argument of lemma \ref{lem1.3}, $length(C_{i})=length(C_{j}) \forall i,j\in \{1,2,\dots,m\}$. Its clear from the definition \ref{df7.1} that there are three cycles of type $H_{1}$ at a vertex. Let $L_{1}(v), L_{2}(v), L_{3}(v)$ be three cycles at $v$. We define $(r, s, k)$-representation of $M$ by the similar idea done in Section 2. Here we take second cut along the cycle $L_{3}$ where the starting adjacent face to base horizontal cycle $L_{1}$ is $6$-gon. Hence we get a $(r, s, k)$-representation of $M$. Next we use the similar argument of theorem \ref{thm4.1}. Here we use the cycle of type $H_{1}$ and $12$-gon and $4$-gon in place of $6$-gon and $3$-gon respectively in the theorem \ref{thm4.1}. Hence we get a set, say $S$, of $12$-gons and $4$-gons. The geometric carrier $|S|$ is a $2$-disk which is bounded by Hamiltonian cycle. Therefore the map of type $\{ 4, 6, 12\}$ on the torus is Hamiltonian.
\hfill$\Box$

\begin{theo} \label{thm7.1} Let $M$ be a map of type $\{4, 6, 12\}$ on the torus. Then it is Hamiltonian.
\end{theo}

\noindent{\sc Proof of Theorem}\ref{thm7.1} The proof of the theorem \ref{thm7.1} follows from the result of this section.
\hfill$\Box$

\section{Map of type $\{ 6, 4, 3, 4\}$}

Let $M$ be a map of type $\{ 6, 4, 3, 4\}$. We define a path as follows $\colon$

\begin{defn}\label{df8.1} Let $P(u_{1}\rightarrow u_{r})$ be a path. The $P$ is said to be well defined if it follow the following properties $\colon$
\begin{enumerate}

\item if $P \cap C_{F_{3}} \neq \emptyset$ then $P\cap C_{F_{3}}$ is a path of length one,

\item if $P \cap C_{F_{4}} \neq \emptyset$ then $P\cap C_{F_{4}}$ is a path of length one,

\item if $P \cap C_{F_{6}} \neq \emptyset$ then $P\cap C_{F_{6}}$ is a path of length two.

\item The boundary vertex of $P$ also well defined in extended path of $P$.

\end{enumerate}
We denote this path by $W_{1}$.
\end{defn}

We define $(r, s, k)$-representation like done in Section 8 of $M$. Here we take second cut along $L_{3}$ where the starting adjacent face to base horizontal cycle $L_{1}$ is $4$-gon. Now we follow the theorem \ref{thm4.1}. Here we consider the cycles of type $W_{1}$ and also $4$-gon in place of $3$-gon. Hence we get a set of $6$-gons and $4$-gons such that the geometric carrier is bounded by a Hamiltonian cycle. Therefore, the map of type $\{ 6, 4, 3, 4\}$ on the torus is Hamiltonian.

\begin{theo} \label{thm8.1} Let $M$ be a map of type $\{6, 4, 3, 4\}$ on the torus. Then it is Hamiltonian.
\end{theo}

\noindent{\sc Proof of Theorem}\ref{thm8.1} The proof of the theorem \ref{thm8.1} follows from the result of this section.
\hfill$\Box$

\section{Acknowledgement}
Work of second author is partially supported by SERB, DST grant No. SR/S4/MS:717/10.

\end{document}